\documentclass[10pt]{article}
\usepackage{pb-diagram}
\usepackage{amsmath, amsfonts}
\usepackage{amssymb}
\usepackage{amscd}

\newtheorem{Df}{Definition}[section]
\newtheorem{Te}[Df]{Theorem}
\newtheorem{Po}[Df]{Proposition}
\newtheorem{Cr}[Df]{Corollary}
\newtheorem{Lm}[Df]{Lemma}
\newtheorem{Ca}[Df]{Claim}
\newtheorem{Cn}[Df]{Conjecture}
\newtheorem{Ex}[Df]{Example}
\newtheorem{Rm}[Df]{Remark}

\newcommand{\Bdf}{\begin{Df}}
\newcommand{\Edf}{\end{Df}}
\newcommand{\Bte}{\begin{Te}}
\newcommand{\Ete}{\end{Te}}
\newcommand{\Bpo}{\begin{Po}}
\newcommand{\Epo}{\end{Po}}
\newcommand{\Bcr}{\begin{Cr}}
\newcommand{\Ecr}{\end{Cr}}
\newcommand{\Blm}{\begin{Lm}}
\newcommand{\Elm}{\end{Lm}}
\newcommand{\Bca}{\begin{Ca}}
\newcommand{\Eca}{\end{Ca}}
\newcommand{\Bcn}{\begin{Cn}}
\newcommand{\Ecn}{\end{Cn}}
\newcommand{\Bex}{\begin{Ex}}
\newcommand{\Eex}{\end{Ex}}
\newcommand{\Brm}{\begin{Rm}}
\newcommand{\Erm}{\end{Rm}}
\newcommand{\Bdm}{{\it Proof.}\ }
\newcommand{\Edm}{\rule{2mm}{2mm}}

\begin{document}

\title{\bf{Koszul calculus}}
\author{Roland Berger, Thierry Lambre and Andrea Solotar
\thanks{\footnotesize This work has been partially supported by the projects  UBACYT 20020130100533BA, PIP-CONICET
	11220150100483CO and MATHAMSUD-REPHOMOL. The third author is a
	research member of CONICET (Argentina).}}
\date{}

\maketitle

\begin{abstract}
We present a calculus which is well-adapted to homogeneous quadratic algebras. We define this calculus on Koszul cohomology 
--resp. homology-- by cup products -- resp. cap products. 
The Koszul homology and cohomology are interpreted in terms of derived categories. If the algebra is not Koszul, then 
Koszul (co)homology
 provides different information than Hochschild (co)homology. As an application of our calculus, the Koszul duality for 
 Koszul cohomology algebras is proved for \emph{any} quadratic algebra, and this duality is extended in some sense to 
 Koszul homology. 
So, the true nature of the Koszul duality theorem is independent of any assumption on the quadratic algebra.
 We compute explicitly this calculus on a non-Koszul example.
\end{abstract} 


\tableofcontents

\section{Introduction}
In this paper, a quadratic algebra is an associative algebra defined by \emph{homogeneous} quadratic relations. Since their definition by Priddy~\cite{priddy:kreso}, 
Koszul algebras form a widely studied class of quadratic algebras~\cite{popo:quad}. In his monograph~\cite{manin:quant}, Manin brings out a 
general approach of quadratic algebras -- not necessarily Koszul -- including the fundamental observation that quadratic algebras form a category 
which should be a relevant framework for a noncommutative analogue of projective algebraic geometry. According to this general approach, 
non-Koszul quadratic algebras deserve more attention.

The goal of this paper is to introduce new homological tools for studying quadratic algebras and to give an application to the Koszul duality. These tools consist in a (co)homology theory, called 
Koszul (co)homology, together with products, called Koszul cup and cap products. They are organized in a calculus, called Koszul calculus. 
If two quadratic algebras are isomorphic in Manin's category \cite{manin:quant}, their Koszul calculi are isomorphic. If the quadratic algebra is 
Koszul, then the Koszul calculus is isomorphic to the Hochschild (co)homology endowed with the usual cup and cap products -- called Hochschild calculus. In 
this introduction, we would like to describe the main features of the Koszul calculus and how they are involved in the course of the paper.

In Section 2, we define the Koszul homology $HK_{\bullet}(A,M)$ of a quadratic algebra $A$ with coefficients in a bimodule $M$ by applying the 
functor $M\otimes_{A^e} -$ to the Koszul complex of $A$, analogously for the Koszul cohomology $HK^{\bullet}(A,M)$. 
If $A$ is Koszul, the Koszul complex is a projective resolution of $A$, so that $HK_{\bullet}(A,M)$ (resp. $HK^{\bullet}(A,M)$) is isomorphic to  Hochschild homology $HH_{\bullet}(A,M)$ (resp. Hochschild cohomology $HH^{\bullet}(A,M)$). Restricting the Koszul calculus to $M=A$, we present 
in Section 9 a non-Koszul quadratic algebra $A$ which is such that $HK_{\bullet}(A)\ncong HH_{\bullet}(A)$ and $HK^{\bullet}(A)\ncong HH^{\bullet}(A)$, showing that $HK_{\bullet}(A)$ and $HK^{\bullet}(A)$ provide new invariants associated to the category of quadratic algebras, 
besides those provided by the Hochschild (co)homology. We prove that the Koszul homology (cohomology) is isomorphic to a Hochschild hyperhomology 
(hypercohomology), showing that this new homology (cohomology) becomes natural in terms of derived categories. 

In Section 3 and Section 4, we introduce the Koszul cup and cap products by restricting the definition of the usual cup and cap products on Koszul 
cochains and chains respectively, providing
 differential graded algebras and differential graded bimodules. These products pass to (co)homology. 

For any unital associative algebra $A$, 
the Hochschild cohomology of $A$ with coefficients in $A$ itself, endowed with the cup product and the Gerstenhaber bracket $[-,-]$, is a Gerstenhaber 
algebra~\cite{gerst:cohom}. We organize the Gerstenhaber algebra structure and the Hochschild homology of $A$, endowed with cap products, in 
a Tamarkin-Tsygan calculus of the kind developed in~\cite{tt:calculus, tl:bvcy}. In the Tamarkin-Tsygan calculus, the Hochschild differential $b$ is 
defined in terms of the multiplication $\mu$ of $A$ and the Gerstenhaber bracket by
\begin{equation} \label{fund1}
b(f)=[\mu, f]
\end{equation}
for any Hochschild cochain $f$. 

It seems difficult to see the Koszul calculus as a Tamarkin-Tsygan calculus because the Gerstenhaber bracket \emph{does not make sense on Koszul cochains}. However, this obstruction can be bypassed by the following formula
\begin{equation} \label{fund2}
b_K(f)=-[e_A, f]_{\underset{K}{\smile}}
\end{equation}
where $b_K$ is the Koszul differential, $e_A$ is the Koszul 1-cocycle defined as the restriction of the Euler derivation $D_A$ of $A$, and $f$ is any Koszul cochain.

In Formula (\ref{fund2}), the symbol $[-, -]_{\underset{K}{\smile}}$ stands for the graded bracket associated to the Koszul cup product $\underset{K}{\smile}$, so that \emph{the Koszul differential may be defined from the Koszul cup product}. 
The Koszul calculus is more flexible than the usual calculus since the formula (\ref{fund2}) is valid for any bimodule $M$, while the definition of the Gerstenhaber bracket is meaningless when considering other bimodules of coefficients~\cite{gs:simplicial}; it is also more symmetric since there 
is an analogue of (\ref{fund2}) in homology, where the Koszul cup product is replaced by the Koszul cap product.

In the Tamarkin-Tsygan calculus, the homology of the Hochschild homology $HH_{\bullet}(A)$ endowed with the Connes differential plays the role of a (generalized) de Rham cohomology of $A$. Since the quadratic algebra $A$ is $\mathbb{N}$-graded and connected, $A$ is acyclic in characteristic zero for this de Rham cohomology (Theorem \ref{poincarelemma}). We give the following Koszul analogue: if $A$ is Koszul, $A$ is acyclic for the \emph{higher Koszul homology}, where we define the higher Koszul homology as the homology of the Koszul homology endowed with the left Koszul cap product by the Koszul class of $e_A$ (Theorem \ref{hkhkoszul}). However, if $A$ is the algebra in the non-Koszul example of Section 9, we prove that $A$ is not acyclic for the higher Koszul homology (Proposition \ref{hkhexample}). Thus the higher Koszul homology is a new invariant of the non-Koszul algebra $A$. We conjecture that the Koszul algebras are exactly the acyclic objects of the higher Koszul homology. 

In~\cite{tl:bvcy}, the second author defined the Tamarkin-Tsygan calculi \emph{with duality}. Specializing this general definition to the 
Hochschild situation, the Tamarkin-Tsygan calculus of an associative algebra $A$ is said to be with duality if there is a class $c$ in a 
space $HH_n(A)$, called the fundamental Hochschild class, such that the $k$-linear map
$$-\frown c: HH^p(A) \longrightarrow HH_{n-p}(A)$$
is an isomorphism for any $p$. If the algebra $A$ is $n$-Calabi-Yau~\cite{vg:cy}, such a calculus exists, and for any bimodule $M$, 
$$-\frown c: HH^p(A,M) \longrightarrow HH_{n-p}(A,M)$$
is an isomorphism coinciding with the Van den Bergh duality~\cite{vdb:dual, tl:bvcy}. Consequently, if $A$ is an $n$-Calabi-Yau Koszul quadratic algebra in characteristic zero, the higher Koszul cohomology 
of $A$ vanishes in any homological degree $p$, except for $p=n$ for which it is one-dimensional (Corollary \ref{hkcohkoszul}). This last fact does not hold for a certain Koszul algebra $A$ of finite global dimension and not Calabi-Yau (Proposition \ref{tensoralgebra}).

In Remark 5.4.10 of~\cite{vg:cy}, Ginzburg mentioned that the Hochschild cohomology algebras of $A$ and its Koszul dual $A^!$ are isomorphic if 
the quadratic algebra $A$ is Koszul. 
This isomorphism of graded algebras was already announced by Buchweitz in the Conference on Representation Theory held in Canberra in 2003, and it has been generalized by Keller in \cite{keller:duality}.
As an application of the Koszul calculus, we obtain such a Koszul duality theorem linking the Koszul cohomology algebras of $A$ and $A^!$ \emph{for any 
quadratic algebra $A$, either Koszul or not} (Theorem \ref{cohkduality}), revealing that the true nature of the Koszul duality theorem is independent of any assumption on quadratic algebras. Our proof of Theorem \ref{cohkduality} uses some standard facts on duality of finite dimensional vector spaces, allowing us to define the Koszul dual of a Koszul cochain (Definition \ref{varphi}).   

Our proof shows two phenomena that already hold for the Koszul algebras. Firstly, the homological weight $p$ is changed by the duality into
 the coefficient weight $m$. Secondly, the 
exchange $p\leftrightarrow m$ implies that we have to replace one of both cohomologies by a 
modified version of the Koszul cohomology and of the Koszul cup product, denoted by tilde accents. 
The statement of Theorem \ref{cohkduality} is the following. 
\Bte \label{cohdual}
Let $V$ be a finite dimensional $k$-vector space and $A=T(V)/(R)$ be a quadratic algebra. Let $A^!=T(V^{\ast})/(R^{\perp})$ be the Koszul dual of 
$A$. There is an isomorphism of $\mathbb{N}\times \mathbb{N}$-graded unital associative algebras
\begin{equation} 
(HK^{\bullet}(A), \underset{K}{\smile}) \cong (\tilde{HK}^{\bullet}(A^!), \tilde{\underset{K}{\smile}}).
\end{equation}
In particular, for any $p\geq 0$ and $m\geq 0$, there is a $k$-linear isomorphism
\begin{equation}
HK^p(A)_m \cong \tilde{HK}^m(A^!)_p.
\end{equation}
\Ete

We illustrate this theorem by an example, with direct computations. Both phenomena are shown to be essential in this example. Theorem \ref{cohkduality} is completed by a bimodule isomorphism in which $HK^{\bullet}(A)$ acts on $HK_{\bullet}(A)$ by cap products (Theorem \ref{hkduality}). 

In Section 9, we compute the Koszul calculus on an example of non-Koszul quadratic algebra $A$. Moreover, we prove that the Koszul 
homology (cohomology) of $A$ is not isomorphic to the Hochschild homology (cohomology) of $A$. For computing the Hochschild homology and cohomology 
of $A$ in degrees 2 and 3, we use a projective bimodule resolution of $A$ due to the third author and Chouhy~\cite{cs:reso}.

\setcounter{equation}{0}

\section{Koszul homology and cohomology}

Throughout the paper, we denote by $k$ the base field and we fix a $k$-vector space $V$. The symbol $\otimes$ will mean $\otimes_k$. The tensor algebra $T(V)=\bigoplus_{m\geq 0} V^{\otimes m}$ of $V$ is graded by the \emph{weight} $m$. For any subspace $R$ of $V^{\otimes 2}$, the associative $k$-algebra $A=T(V)/(R)$ is called a \emph{quadratic algebra}, and it inherits the grading by the weight.
We denote the homogeneous component of weight $m$ of $A$ by $A_m$.

\subsection{Recalling the bimodule complex $K(A)$}

Let $A=T(V)/(R)$ be a quadratic algebra. For the definition of the bimodule complex $K(A)$, we follow Van den Bergh, precisely Section 3 of~\cite{vdb:hom}. Notice that our $K(A)$ is denoted by $K'(A)$ in~\cite{vdb:hom}. For any $p\geq 2$, we define the subspace  $W_p$ of $V^{\otimes p}$ by 
$$W_{p}=\bigcap_{i+2+j=p}V^{\otimes i}\otimes R\otimes V^{\otimes j}, \hbox{ where } i, j\ge 0,$$

\noindent while $W_0=k$ and $W_1=V$. It is convenient to use the following notation: an 
arbitrary element of $W_p$ will be denoted by a product $x_1 \ldots  x_p$, where $x_1, \ldots , x_p$ are in $V$. This notation should be thought of 
as a sum of such products. Moreover, regarding $W_p$ as a subspace of $V^{\otimes q}\otimes W_r \otimes V^{\otimes s}$ where $q+r+s=p$, the 
element $x_1 \ldots  x_p$ viewed in $V^{\otimes q}\otimes W_r \otimes V^{\otimes s}$ will be denoted by the \emph{same} notation, meaning that 
$x_{q+1} \ldots x_{q+r}$ is thought of as a sum belonging to $W_r$ and the other $x$'s are thought of as arbitrary elements in $V$. \emph{We will systematically use this notation throughout the paper.}

Clearly, $V$ is the component of weight 1 of $A$, so that $V^{\otimes p}$ is a subspace of $A^{\otimes p}$. As defined by 
Van den Bergh~\cite{vdb:hom}, the \emph{Koszul complex} $K(A)$ of the quadratic algebra $A$ is a weight graded bimodule subcomplex of the bar resolution
$B(A)$ of $A$. Precisely, $K(A)_p=K_p$ is the subspace $A\otimes W_p \otimes A$ of 
$A\otimes A^{\otimes p} \otimes A$. It is easy to see that $K(A)$ coincides with the complex
\begin{equation} \label{priddy}
\cdots \stackrel{d}{\longrightarrow} K_{p} \stackrel{d}{\longrightarrow} K_{p-1} \stackrel{d}{\longrightarrow} \cdots
\stackrel{d}{\longrightarrow} K_{1} \stackrel{d}{\longrightarrow} K_{0} \stackrel{}{\longrightarrow} 0\,,
\end{equation}
where the differential $d$ is defined on $K_p$ as follows
\begin{equation} \label{defd}
d(a \otimes x_1 \ldots x_{p} \otimes a') =ax_1\otimes x_2 \ldots x_{p}\otimes a'+(-1)^p a \otimes x_1 \ldots x_{p-1}\otimes x_{p} a',
\end{equation}
for $a$, $a'$ in $A$ and $x_1 \ldots x_{p}$ in $W_p$, using the above notation. 
The homology of $K(A)$ is equal to $A$ in degree $0$, and to $0$ in degree $1$. The following definition takes into account the bimodule complex $K(A)$, instead of the left or right module versions of the Koszul complex commonly used for defining Koszul algebras~\cite{popo:quad,lv:alop}. The following definition is equivalent to the usual one, according to Proposition 3.1 in~\cite{vdb:hom} and its obvious converse. 
\Bdf \label{defK}
A quadratic algebra $A$ is said to be Koszul if the homology of $K(A)$ is $0$ in any positive degree.
\Edf

The multiplication $\mu: K_0 =A\otimes A \rightarrow A$ defines a morphism from the complex $K(A)$ to the complex $A$ concentrated in degree $0$.
Whereas $\mu: B(A) \rightarrow A$ is always a quasi-isomorphism,
$A$ is Koszul if and only if $\mu: K(A)\rightarrow A$ is a quasi-isomorphism. So, if the quadratic algebra $A$ is Koszul, the bimodule 
free resolution $K(A)$ may be used to compute the Hochschild homology and cohomology of $A$ instead of $B(A)$. In the two subsequent subsections, the same (co)homological functor is defined by replacing $B(A)$ by $K(A)$ even if $A$ is not Koszul. The goal of this paper is to show that the so-obtained Koszul (co)homology 
is of interest for quadratic algebras, providing invariants that are not obtained with the Hochschild (co)homology.

\subsection{The Koszul homology $HK_{\bullet}(A,M)$} \label{koh}

Let $M$ be an $A$-bimodule. As usual, $M$ can be considered as a left or right $A^e$-module, where $A^e=A\otimes A^{op}$. Applying the 
functor $M\otimes_{A^e} -$ to $K(A)$, we obtain the chain complex ($M\otimes W_{\bullet}$, $b_K$), where $W_{\bullet}=\bigoplus_{p\geq 0} W_p$. 
The elements of $M\otimes W_p$ are called the \emph{Koszul $p$-chains with coefficients in $M$}. From Equation (\ref{defd}), we see that 
the differential $b_K=M\otimes_{A^e} d$ is given on $M\otimes W_p$ by the formula
\begin{equation} \label{defb}
b_K(m \otimes x_1 \ldots x_{p}) =m.x_1\otimes x_2 \ldots x_{p} +(-1)^p x_p .m \otimes x_1 \ldots x_{p-1},
\end{equation}
for any $m$ in $M$ and $x_1 \ldots x_{p}$ in $W_p$, using the notation of Subsection 2.1.
\Bdf \label{hko}
Let $A=T(V)/(R)$ be a quadratic algebra and $M$ be an $A$-bimodule. The homology of the complex ($M\otimes W_{\bullet}$, $b_K$) is called the 
Koszul homology of $A$ with coefficients in $M$, and is denoted by $HK_{\bullet}(A,M)$. We set $HK_{\bullet}(A)=HK_{\bullet}(A,A)$.
\Edf

The inclusion  $\chi:K(A)\rightarrow B(A)$ induces a morphism of complexes
$\tilde{\chi}=M\otimes_{A^e}\! \chi$ from $(M\otimes W_{\bullet}, b_K)$ to 
$(M\otimes A^{\otimes \bullet}, b)$, where $b$ is the Hochschild differential. For each degree $p$, $\tilde{\chi}_p$ coincides with the 
natural injection of $M\otimes W_p$ into $M\otimes A^{\otimes p}$. Since the complex
\begin{equation} \label{debd}
A\otimes R\otimes A \stackrel{d}\longrightarrow
A\otimes  V\otimes A \stackrel{d}\longrightarrow A\otimes  A \stackrel{\mu}\longrightarrow A \rightarrow 0
\end{equation} is exact, the $k$-linear map $H(\tilde{\chi})_p:HK_p(A,M)\rightarrow HH_p(A,M)$ is an isomorphism for $p=0$ and $p=1$. The following is clear.
\Bpo  \label{hkhh}
Let $A=T(V)/(R)$ be a Koszul quadratic algebra. For any $A$-bimodule $M$ and any $p\geq 0$, 
$H(\tilde{\chi})_p$ is an isomorphism.  \Epo

For the non-Koszul algebra $A$ of Section 9, we will see that $H(\tilde{\chi})_3$ is not surjective when $M=A$. Quadratic $k$-algebras form Manin's category~\cite{manin:quant}. In this category, a morphism 
$u$ from $A=T(V)/(R)$ to $A'=T(V')/(R')$ is determined by a linear map $u : V \rightarrow V'$ such that $u ^{\otimes 2} (R) \subseteq R'$.
For each $p$, $u ^{\otimes p}$ maps $W_p$ into $W'_p$, with obvious notation. Moreover, the maps 
$a\otimes  x_1 \ldots x_{p} \mapsto u (a) \otimes  u(x_1) \ldots u(x_{p})$ define a morphism of complexes from 
$(A\otimes W_{\bullet}, b_K)$ to $(A'\otimes W'_{\bullet}, b_K)$. So we obtain a covariant functor $A \mapsto HK_{\bullet}(A)$.\\

Let us now show that the Koszul homology is isomorphic to a Hochschild hyperhomology, namely
\begin{equation} \label{hyperhh}
HK_{\bullet}(A,M) \cong \mathbb{H}\mathbb{H}_{\bullet}(A, M\otimes _A K(A)).
\end{equation}
Denote by $\mathcal{A}$ (resp. $\mathcal{E}$) the abelian category of $A$-bimodules (resp. $k$-vector spaces). For any $A$-bimodule $M$, the 
left derived functor $M \stackrel{L}\otimes_{A^e}-$ is defined from the triangulated category $\mathcal{D}^- (\mathcal{A})$ to the 
triangulated category $\mathcal{D}^- (\mathcal{E})$, so that we have
\begin{equation} \label{dercat}
HK_{p}(A,M) \cong H_p (M \stackrel{L}\otimes_{A^e} K(A)).
\end{equation}
The following lemma is standard, used e.g. in the proof of the Van den Bergh duality~\cite{vdb:dual}. 
\Blm  \label{hlem}
Let $M$ and $N$ be $A$-bimodules. The $k$-linear map
$$\zeta: M\otimes_{A^e} N \rightarrow (M\otimes _A N)\otimes_{A^e} A$$
defined by $\zeta (x\otimes_{A^e} y)=(x\otimes _A y)\otimes_{A^e} 1$
is an isomorphism. Moreover, for any complex of $A$-bimodules $C$, the map $\zeta: M\otimes_{A^e} C \rightarrow (M\otimes _A C)\otimes_{A^e} A$ 
is an isomorphism of complexes.
\Elm

\medskip

In other words, the functor $F: C\mapsto M\otimes_{A^e} C$ coincides with the composite $H\circ G$ where $G: C\mapsto M\otimes _A C$ and 
$H: C' \mapsto C'\otimes_{A^e} A$. So their left derived functors satisfy $LF\cong L(H)\circ L(G)$, in particular for $C=K(A)$,
\begin{equation} \label{isostand}
M \stackrel{L}\otimes_{A^e} K(A) \cong (M \stackrel{L}\otimes_{A} K(A)) \stackrel{L}\otimes_{A^e} A.
\end{equation}
Passing to homology and using the definition of hypertor~\cite{weib:homo}, we obtain
\begin{equation} \label{hypertor}
HK_p(A,M) \cong \mathbb{T}or_p^{A^e}(M\otimes _A K(A), A),
\end{equation}
which proves the isomorphism (\ref{hyperhh}). If $A$ is Koszul, we recover usual $Tor$ and Proposition \ref{hkhh}.

\subsection{The Koszul cohomology $HK^{\bullet}(A,M)$} \label{kocoh}

Throughout, $Hom_k$ will be denoted by $Hom$. Applying the functor $Hom_{A^e} (-, M)$ to the complex $K(A)$, 
we obtain the cochain complex ($Hom(W_{\bullet},M)$, $b_K$), where
$$Hom(W_{\bullet},M)=\bigoplus_{p\geq 0}Hom(W_p,M).$$
The elements of $Hom(W_p,M)$ are called the \emph{Koszul $p$-cochains with coefficients in $M$}. Given a Koszul $p$-cochain $f:W_p\rightarrow M$, 
its differential $b_K(f)=-(-1)^p f\circ d$ is defined by 
\begin{equation} \label{defcob}
b_K(f)( x_1 \ldots x_{p+1}) =f(x_1\ldots x_{p}).x_{p+1} -(-1)^p x_1 .f(x_2 \ldots x_{p+1}),
\end{equation}
for any $x_1 \ldots x_{p+1}$ in $W_{p+1}$, using the notation of Subsection 2.1. 
\Bdf \label{cohko}
Let $A=T(V)/(R)$ be a quadratic algebra and $M$ an $A$-bimodule. The homology of the complex ($Hom(W_{\bullet},M)$, $b_K$) is called the 
Koszul cohomology of $A$ with coefficients in $M$, and is denoted by $HK^{\bullet}(A,M)$. We set $HK^{\bullet}(A)=HK^{\bullet}(A,A)$.
\Edf

The map $\chi^{\ast}=Hom_{A^e} (\chi, M)$ defines a morphism 
of complexes from $(Hom(A^{\otimes \bullet},M), b)$ to $(Hom(W_{\bullet},M), b_K)$, where $b$ is the Hochschild differential. For each degree $p$, $\chi^{\ast}_p$ coincides with the 
 natural projection of $Hom(A^{\otimes p},M)$ onto $Hom(W_p,M)$. The $k$-linear 
 map $H(\chi^{\ast})_p:HH^p(A,M)\rightarrow HK^p(A,M)$ is an isomorphism for $p=0$ and $p=1$. 
\Bpo  \label{cohkhh}
Let $A=T(V)/(R)$ be a Koszul quadratic algebra. For any $A$-bimodule $M$ and any $p\geq 0$, 
$H(\chi^{\ast})_p$ is an isomorphism.  \Epo

In the non-Koszul example of Section 9, we will see that $H(\chi^{\ast})_2$ is not surjective for $M=A$. Here again, the same functorial properties of Hochschild cohomology stand for 
Koszul cohomology. In particular, there is a contravariant functor $A \mapsto HK^{\bullet}(A,A^{\ast})$, where the $A$-bimodule $A^{\ast}=Hom(A,k)$ 
is defined by: $(a.f.a')(x)=f(a'xa)$ for any $k$-linear map $f:A\rightarrow k$, and $x$, $a$, $a'$ in $A$.

As we prove now, the Koszul cohomology is isomorphic to the following Hochschild hypercohomology
\begin{equation} \label{hyperhcoh}
HK^{\bullet}(A,M) \cong \mathbb{H}\mathbb{H}^{\bullet}(A, Hom _A (K(A),M)).
\end{equation}
For any $A$-bimodule $M$, the right derived functor $RHom_{A^e}(-, M)$ is defined from the triangulated category $\mathcal{D}^- (\mathcal{A})$ to 
the triangulated category $\mathcal{D}^+ (\mathcal{E})$, so that we have
\begin{equation} \label{codercat}
HK^{p}(A,M) \cong H^p (RHom_{A^e}(K(A), M)).
\end{equation}
The proof continues as in homology by using the next lemma. We leave details to the reader.
\Blm  \label{cohlm}
Let $M$ and $N$ be $A$-bimodules. The $k$-linear map
$$\eta: Hom_{A^e}(N, M) \rightarrow Hom_{A^e} (A,Hom_A(N,M))$$
defined by $\eta (f)(a)(x)=f(xa)$ for any $A$-bimodule map $f:N\rightarrow M$, $a$ in $A$ and $x$ in $N$,
is an isomorphism, where $Hom_A(N,M)$ denotes the space of left $A$-module morphisms from $M$ to $N$. Moreover, for any complex of $A$-bimodules $C$, 
$\eta: Hom_{A^e} (C,M) \rightarrow Hom_{A^e} (A, Hom_A(C,M))$ is an isomorphism of complexes.
\Elm

\subsection{Coefficients in $k$} \label{const}

In this subsection, the Koszul homology and cohomology are examined for the trivial bimodule $M=k$. Denote by $\epsilon : A \rightarrow k$ the 
augmentation of $A$, so that the $A$-bimodule $k$ is defined by the following actions: $a.1.a'=\epsilon (aa')$ for 
any $a$ and $a'$ in $A$. It is immediate from (\ref{defb}) and (\ref{defcob}) that the differentials $b_K$ vanish in case $M=k$. Denoting $Hom(E,k)$ 
by $E^{\ast}$ for any $k$-vector space $E$, we obtain the following.
\Bpo \label{hkk}
Let $A=T(V)/(R)$ be a quadratic algebra. For any $p\geq 0$, we have $HK_p(A,k) = W_p$ and $HK^p(A,k) = W^{\ast}_p$.
\Epo

Let us give a conceptual explanation of this proposition. We consider quadratic algebras as connected algebras graded by the weight~\cite{popo:quad}. 
Let $A=T(V)/(R)$ be a quadratic algebra. In the category of graded $A$-bimodules, $A$ has a minimal projective resolution $P(A)$ whose component 
of homological degree $p$ has the form $A\otimes E_p \otimes A$, where $E_p$ is a weight graded space. Moreover, the minimal weight in $E_p$ is 
equal to $p$ and the component of weight $p$ in $E_p$ coincides with $W_p$. Denote by $\underline{Hom}$ the graded $Hom$ w.r.t. the weight grading 
of $A$, and by $\underline{HH}$ the corresponding graded Hochschild cohomology. The following fundamental property of $P(A)$ holds for any 
connected graded algebra $A$.
\Blm \label{hhk}
The differentials of the complexes $k\otimes_{A^e} P(A)$ and $\underline{Hom}_{A^e} (P(A), k)$ vanish.
\Elm

Consequently, there are isomorphisms $HH_p(A,k) \cong E_p$ and $\underline{HH}^p(A,k) \cong \underline{Hom}(E_p, k)$ for any $p\geq 0$. Since $K(A)$ is 
a weight graded subcomplex of $P(A)$, $H(\tilde{\chi})_p$ coincides with the natural injection of $W_p$ into $E_p$ and $H(\chi^{\ast})_p$ with the natural projection of $\underline{Hom}(E_p, k)$ onto $W^{\ast}_p$. So we obtain the following converses of Proposition \ref{hkhh} and Proposition \ref{cohkhh}.
\Bpo \label{conv}
Let $A=T(V)/(R)$ be a quadratic algebra. The algebra $A$ is Koszul if either (i) or (ii) hold.

(i) For any $p\geq 0$, $H(\tilde{\chi})_p:HK_p(A,k)\rightarrow HH_p(A,k)$ is an isomorphism.

(ii) For any $p\geq 0$, $H(\chi^{\ast})_p:\underline{HH}^p(A,k)\rightarrow HK^p(A,k)$ is an isomorphism.
\Epo

\setcounter{equation}{0}

\section{The Koszul cup product}

\subsection{Definition and first properties}

We define the Koszul cup product $\underset{K}{\smile}$ of Koszul cochains by restricting the usual cup product $\smile$ of Hochschild cochains recalled e.g. in~\cite{tl:bvcy}. We use the notation of Subsection 2.1.
\Bdf \label{defcup}
Let $A=T(V)/(R)$ be a quadratic algebra. Let $P$ and $Q$ be $A$-bimodules. For any Koszul $p$-cochain $f:W_p\rightarrow P$ and any Koszul $q$-cochain $g:W_q\rightarrow Q$, define the 
Koszul $(p+q)$-cochain $f\underset{K}{\smile} g : W_{p+q} \rightarrow P\otimes_A Q$ by the following equality
\begin{equation} \label{kocup}
(f\underset{K}{\smile} g) (x_1 \ldots x_{p+q}) = (-1)^{pq} f(x_1 \ldots x_p)\otimes_A \, g(x_{p+1} \ldots  x_{p+q}),
\end{equation}
for any $x_1 \ldots x_{p+q} \in W_{p+q}$.
\Edf

The Koszul cup product $\underset{K}{\smile}$ is $k$-bilinear and associative, and we have the formula
\begin{equation} \label{hokocup}
\chi^{\ast}(F \smile G) = \chi^{\ast}(F) \underset{K}{\smile} \chi^{\ast}(G)
\end{equation}
for any Hochschild cochains $F:A^{\otimes p} \rightarrow P$ and $G:A^{\otimes q} \rightarrow Q$. We deduce the identity
\begin{equation} \label{kodga}
b_K(f \underset{K}{\smile} g)= b_K(f) \underset{K}{\smile} g +(-1)^p f \underset{K}{\smile} b_K(g),
\end{equation}
from the identity known for the usual $\smile$. In particular, $Hom(W_{\bullet},A)$ is a differential graded algebra (dga). For any $A$-bimodule $M$, $Hom(W_{\bullet},M)$ is a differential graded bimodule over the dga $Hom(W_{\bullet},A)$. 
The proof of the following statement is clear.
\Bpo
Let $A=T(V)/(R)$ be a quadratic algebra. The Koszul cup product $\underset{K}{\smile}$ defines a Koszul cup product, still denoted by $\underset{K}{\smile}$, on Koszul cohomology 
classes. A formula similar to (\ref{hokocup}) holds for $H(\chi^{\ast})$. 
Endowed with this product, $HK^{\bullet}(A)$ and $HK^{\bullet}(A,k)$ are graded associative algebras. For any $A$-bimodule $M$, $HK^{\bullet}(A,M)$ is a graded $HK^{\bullet}(A)$-bimodule.
\Epo

Since $HK^0(A)=Z(A)$ is the center of the algebra $A$, $HK^{\bullet}(A,M)$ is a $Z(A)$-bimodule. From Proposition \ref{hkk}, $HK^{\bullet}(A,k)$ coincides with the graded 
algebra $W_{\bullet}^{\ast}=\bigoplus_{p\geq 0} W_p^{\ast}$ endowed with the graded tensor product of linear forms composed with inclusions $W_{p+q}\hookrightarrow W_p\otimes W_q$. 
Recall that the graded algebra $(\underline{HH}^{\bullet}(A,k),\smile)$ is isomorphic to the Yoneda algebra $E(A)=\underline{Ext}_A^{\ast}(k,k)$ of the graded algebra 
$A$~\cite{popo:quad}. 
\Bpo
Let $A=T(V)/(R)$ be a quadratic algebra. The map $H(\chi^{\ast})$ defines a graded algebra morphism from the Yoneda algebra $E(A)$ of $A$ onto $W_{\bullet}^{\ast}$, and this is an isomorphism 
if and only if $A$ is Koszul.
\Epo

\subsection{The Koszul cup bracket}

\Bdf \label{defcupbra}
Let $A=T(V)/(R)$ be a quadratic algebra. Let $P$ and $Q$ be $A$-bimodules, at least one of them equal to $A$. For any Koszul $p$-cochain $f:W_p\rightarrow P$ and any Koszul 
$q$-cochain $g:W_q\rightarrow Q$, we define the Koszul cup bracket 
by
\begin{equation} \label{kocupbra}
[f, g]_{\underset{K}{\smile}} =f\underset{K}{\smile} g - (-1)^{pq} g\underset{K}{\smile} f.
\end{equation}
\Edf

The Koszul cup bracket is $k$-bilinear, graded antisymmetric, and it passes to cohomology. 
We still use the notation $[\alpha, \beta]_{\underset{K}{\smile}}$ for the cohomology classes $\alpha$ and $\beta$ of $f$ and $g$. The Koszul cup bracket is a graded biderivation 
of the graded associative algebras $Hom(W_{\bullet},A)$ and $HK^{\bullet}(A)$. We will see that the Koszul cup bracket plays in some sense the role of the Gerstenhaber bracket. For this,
we will consider the Euler derivation of $A$ as a Koszul 1-cocycle.

\subsection{The fundamental 1-cocycle}

\Blm
Let $A=T(V)/(R)$ be a quadratic algebra. Let $f: V \rightarrow V$ be a $k$-linear map considered as a Koszul 1-cochain with coefficients in $A$. If $f$ is 
a coboundary, then $f=0$. If $f$ is a cocycle, then its cohomology class contains a unique $1$-cocycle with image in $V$ and this cocycle is equal to $f$. \Elm
\Bdm
If $f=b_K(a)$ for some $a$ in $A$, then $f(x)=ax-xa$ for any $x$ in $V$. Since $f(x) \in V$, this implies that $f(x)=a_0x-xa_0$ with $a_0\in k$, thus $f=0$.  
\Edm

\Bdf
Let $A=T(V)/(R)$ be a quadratic algebra. The Euler derivation --also called weight map-- $D_A: A \rightarrow A$ of the graded algebra $A$ is defined by
$D_A(a)=ma$
for any $m\geq 0$ and any homogeneous element $a$ of weight $m$ in $A$.  
\Edf

We denote by $e_A$ the restriction of $D_A$ to $V$. The map $e_A: V \rightarrow A$ is a Koszul 1-cocycle called the \emph{fundamental $1$-cocycle} of $A$. It is defined by $e_A(x)=x$ for any $x$ in $V$. It corresponds to the canonical element $\xi_A$ of Manin~\cite{manin:quant}. By the previous lemma, $e_A$ is not a coboundary if $V \neq 0$. The Koszul class of $e_A$ is denoted by $\overline{e}_A$ 
and it is called the \emph{fundamental $1$-class} of $A$. The following statement is easily proved, but it is of crucial importance for the Koszul calculus.

\Bte \label{thfundacoho}
Let $A=T(V)/(R)$ be a quadratic algebra. For any Koszul cochain $f$ with coefficients in any $A$-bimodule $M$, the following formula holds
\begin{equation} \label{fundacoho}
[e_A, f]_{\underset{K}{\smile}}=-b_K(f).
\end{equation}
\Ete
\Bdm
For any $x_1 \ldots x_{p+1}$ in $W_{p+1}$, one has $(e_A\underset{K}{\smile} f) (x_1 \ldots x_{p+1}) = (-1)^p x_1 . f(x_2 \ldots  x_{p+1})$ and 
$(f\underset{K}{\smile} e_A) (x_1 \ldots x_{p+1}) = (-1)^p f(x_1 \ldots x_p).x_{p+1}$, so that Formula (\ref{fundacoho}) is immediate from (\ref{defcob}). 
\Edm
\\

The fundamental formula (\ref{fundacoho}) shows that \emph{the Koszul differential $b_K$ may be defined from the Koszul cup product}, and doing so, we may deduce the identity (\ref{kodga}) from the biderivation $[-, -]_{\underset{K}{\smile}}$. The simple formula (\ref{fundacoho}) is 
replaced in the Hochschild calculus by the ``more sophisticated'' and well-known formula
\begin{equation} \label{fundagerst}
b(F)=[\mu, F],
\end{equation}
where $[-, -]$ is the Gerstenhaber bracket, multiplication $\mu=b(Id_A)$ is a $2$-coboundary and $F$ is any Hochschild cochain.

Let us show that it is possible to deduce the fundamental formula (\ref{fundacoho}) from the Gerstenhaber calculus, that is, from the Hochschild calculus including the Gerstenhaber 
product $\circ$. We recall from~\cite{gerst:cohom} the Gerstenhaber identity  
\begin{equation} \label{gerstidentity}
b(F\circ G)=b(F)\circ G - (-1)^p F\circ b(G) - (-1)^p [F, G]_{\smile}
\end{equation}
for any Hochschild cochains $F:A^{\otimes p} \rightarrow A$ and $G:A^{\otimes q} \rightarrow A$, where
$$[F, G]_{\smile} =F \smile G - (-1)^{pq} G \smile F.$$
The Gerstenhaber product $F\circ G$ is the $(p+q-1)$-cochain defined by
\small
\begin{equation} \label{gerstproduct}
F\circ G\,(a_1, \ldots, a_{p+q-1})= \sum_{1\leq i \leq p} (-1 )^{(i-1)(q-1)} F(a_1, \ldots a_{i-1}, G(a_i, \ldots a_{i+q-1}), a_{i+q}, \ldots , a_{p+q-1}),
\end{equation}
\normalsize
for any $a_1, \ldots, a_{p+q-1}$ in $A$.

For $G=D_A$, identity (\ref{gerstidentity}) becomes
$$b(F\circ D_A)-b(F)\circ D_A = - (-1)^p [F, D_A]_{\smile}.$$
Restricting this identity to $W_{p+1}$, the right-hand side coincides with $[e_A, f]_{\underset{K}{\smile}}$ where $f$ is the restriction of $F$ to $W_p$. Since $F\circ D_A= pf$ on $W_p$, 
the restriction of $b(F\circ D_A)$ is equal to $pb_K(f)$. The restriction of $b(F)\circ D_A$ is equal to $(p+1)b_K(f)$. Thus we recover the fundamental formula $[e_A, f]_{\underset{K}{\smile}}=-b_K(f)$.

\subsection{Koszul derivations} \label{subseckoder}

\Bdf \label{defkoder}
Let $A=T(V)/(R)$ be a quadratic algebra and let $M$ be an $A$-bimodule. Any Koszul 1-cocycle $f: V \rightarrow M$ with coefficients in $M$ will be called a Koszul derivation of $A$ 
with coefficients in $M$. When $M=A$, we will simply speak about a Koszul derivation of $A$. \Edf

According to equation (\ref{defcob}), a $k$-linear map $f: V \rightarrow M$ is a Koszul derivation if and only if
\begin{equation} \label{koder}
f(x_1)x_2 +x_1f(x_2)=0,
\end{equation}
for any $x_1x_2$ in $R$ (using the notation of Subsection 2.1). If this equality holds, the unique derivation $\tilde{f}: T(V) \rightarrow M$ extending $f$ defines a unique derivation 
$D_f: A \rightarrow M$ from the algebra $A$ to the bimodule $M$. The $k$-linear map $f\mapsto D_f$ is an isomorphism from the space of Koszul derivations of $A$ with coefficients in $M$ 
to the space of derivations from $A$ to $M$. As for (\ref{fundacoho}), it is possible to deduce the following from the Gerstenhaber calculus.
\Bpo
Let $A=T(V)/(R)$ be a quadratic algebra and let $M$ be an $A$-bimodule. For any Koszul derivation $f: V \rightarrow M$ and any Koszul $q$-cocycle $g:W_q \rightarrow A$, one has 
\begin{equation} \label{kocupbraspe}
[f, g]_{\underset{K}{\smile}} =b_K(D_f\circ g).
\end{equation}
\Epo
\Bdm
Applying $D_f$ to equation $g(x_1\ldots x_{q}).x_{q+1} =(-1)^q x_1 .g(x_2 \ldots x_{q+1})$, we get
\small
$$D_f(g(x_1\ldots x_{q})).x_{q+1} + g(x_1\ldots x_{q}).f(x_{q+1})=(-1)^q (f(x_1).g(x_2 \ldots x_{q+1})+ x_1.D_f(g(x_2 \ldots x_{q+1}))),$$
\normalsize
and equality (\ref{kocupbraspe}) follows from (\ref{defcob}).
\Edm
\Bcr \label{kocupcommutation}
Let $A=T(V)/(R)$ be a quadratic algebra and let $M$ be an $A$-bimodule. For any $\alpha \in HK^p(A,M)$ with $p=0$ or $p=1$ and  $\beta \in HK^q(A)$, one has the identity
\begin{equation} \label{kocupbrazero}
[\alpha, \beta]_{\underset{K}{\smile}} =0.
\end{equation}
\Ecr
\Bdm
The case $p=1$ follows from the proposition. The case $p=0$ is clear since $HK^0(A,M)$ is the space of the elements of $M$ commuting to any element of $A$.
\Edm
\\

If $A$ is Koszul, then $[\alpha, \beta]_{\underset{K}{\smile}}=0$ for any $p$ and $q$, using the Gerstenhaber calculus and the isomorphisms $H(\chi^{\ast})$. We do not know whether 
$[\alpha, \beta]_{\underset{K}{\smile}}=0$ holds for any $p$ and $q$ when $A$ is not Koszul. It holds for $M=A$ by direct verifications in the non-Koszul example of Section 9. Observe that, in this example, $H(\chi^{\ast})_2$ is not surjective for $M=A$, so that there exists a Koszul $2$-cocycle which does 
not extend to a Hochschild $2$-cocycle. Consequently, it seems hard to prove the identity (\ref{kocupbrazero}) for $p=q=2$ \emph{in general} from the Gerstenhaber calculus. Notice also that the equality (\ref{gerstproduct}) defining the 
Gerstenhaber product \emph{does not make sense} for $f\circ g:W_{p+q-1} \rightarrow A$ when $f:W_p\rightarrow A$ and $g:W_q \rightarrow A$.

\subsection{Higher Koszul cohomology} \label{subsechkoco}

Let $A=T(V)/(R)$ be a quadratic algebra. Let $f: V \rightarrow A$ be a Koszul derivation of $A$. 
Denote by $[f]$ the cohomology class of $f$. Assuming $\mathrm{char}(k)\neq 2$, identity (\ref{kocupbrazero}) shows that $[f]\underset{K}{\smile}[f]=0$, so that the $k$-linear 
map $[f]\underset{K}{\smile} -$ is a cochain differential on $HK^{\bullet}(A,M)$ for any $A$-bimodule $M$. We obtain therefore a new cohomology, called \emph{higher 
Koszul cohomology associated to $f$}. The Gerstenhaber identity (\ref{gerstidentity}) implies that $2 D_f \smile D_f = b(D_f \circ D_f)$, therefore $[D_f ] \smile -$ is 
a cochain differential on $HH^{\bullet}(A,M)$, defining a \emph{higher Hochschild cohomology associated to $f$}. Moreover $H(\chi^{\ast})$ induces a morphism from the higher 
Hochschild cohomology to the higher Koszul cohomology, which is an isomorphism if $A$ is Koszul.

Let us limit ourselves to the case $f=e_A$, the fundamental $1$-cocycle. In this case, without any assumption on the characteristic of $k$, the formula $e_A\underset{K}{\smile}e_A=0$ shows 
that the $k$-linear map $e_A\underset{K}{\smile} -$ is a cochain differential on $Hom(W_{\bullet},M)$, and $\overline{e}_A\underset{K}{\smile} -$ is a cochain differential 
on $HK^{\bullet}(A,M)$.

\Bdf \label{hikoco}
Let $A=T(V)/(R)$ be a quadratic algebra and let $M$ be an $A$-bimodule. The differential $\overline{e}_A\underset{K}{\smile} -$ of $HK^{\bullet}(A,M)$ is denoted by $\partial_{\smile}$. 
The homology of $HK^{\bullet}(A,M)$ endowed with $\partial_{\smile}$ is called the higher Koszul cohomology of $A$ with coefficients in $M$ and is denoted by $HK_{hi}^{\bullet}(A,M)$. 
We set $HK_{hi}^{\bullet}(A)=HK_{hi}^{\bullet}(A,A)$.
\Edf

If we want to evaluate $\partial_{\smile}$ on classes, it suffices to go back to the formula
$$(e_A\underset{K}{\smile} f) (x_1 \ldots x_{p+1}) = f(x_1 \ldots x_p). x_{p+1}$$
for any cocycle $f:W_p \rightarrow M$, and any $x_1 \ldots x_{p+1}$ in $W_{p+1}$. Since $HK^0(A,M)$ equals the space $Z(M)$ of the elements of $M$ commuting to any element of $A$, we obtain the following.
\Bpo  \label{zerohikoco}
Let $A=T(V)/(R)$ be a quadratic algebra and let $M$ be an $A$-bimodule. $HK_{hi}^0(A,M)$ is the space of the elements $u$ of $Z(M)$ such that there exists $v\in M$ satisfying 
$u.x=v.x-x.v$ for any $x$ in $V$. In particular, if the bimodule $M$ is symmetric,
then $HK_{hi}^0(A,M)$ is the space of elements of $M$ annihilated by $V$. If $A$ is a commutative domain and $V\neq 0$, then $HK_{hi}^0(A)=0$.
\Epo

The differential $e_A\underset{K}{\smile} -$ vanishes for $M=k$, hence Proposition \ref{hkk} implies that $HK_{hi}^p(A,k)=W^{\ast}_p$ for any $p\geq 0$.

 \subsection{Higher Koszul cohomology with coefficients in $A$} \label{subsechkcoAcoeff}

\Blm
Let $A=T(V)/(R)$ be a quadratic algebra. Given $\alpha$ in $HK^p(A)$ and $\beta$ in $HK^q(A)$, 
$$\partial_{\smile}(\alpha \underset{K}{\smile} \beta)= \partial_{\smile}(\alpha) \underset{K}{\smile} \beta = (-1)^p \alpha \underset{K}{\smile} \partial_{\smile}(\beta).$$
\Elm
\Bdm
The first equality comes from $\overline{e}_A\underset{K}{\smile}(\alpha \underset{K}{\smile} \beta)=(\overline{e}_A\underset{K}{\smile}\alpha) \underset{K}{\smile} \beta$. The second one is clear from the relation $[\overline{e}_A, \alpha]_{\underset{K}{\smile}}=0$.
\Edm
\\

Consequently, the Koszul cup product is defined on $HK_{hi}^{\bullet}(A)$, still denoted by $\underset{K}{\smile}$, and $(HK_{hi}^{\bullet}(A), \underset{K}{\smile})$ is a 
graded associative algebra. Remark that, if $V \neq 0$, then $\partial_{\smile}(1)=\overline{e}_A \neq 0$ and so $1$ and $\overline{e}_A$ \emph{do not survive} in higher Koszul cohomology. 
To go further in the structure of $HK_{hi}^{\bullet}(A)$, we require a finiteness assumption.

Throughout the remainder of this subsection, assume that $V$ is \emph{finite dimensional}. A Koszul $p$-cochain $f:W_p \rightarrow A_m$ is said to be homogeneous of weight $m$. The space of Koszul cochains $Hom(W_{\bullet},A)$ is 
$\mathbb{N}\times \mathbb{N}$-graded by the \emph{biweight} $(p,m)$, where $p$ is called the \emph{homological weight} and $m$ is called the \emph{coefficient weight}. If $f:W_p \rightarrow A_m$ 
and $g:W_q \rightarrow A_n$ are homogeneous of biweights $(p,m)$ and $(q,n)$ respectively, then $f\underset{K}{\smile} g : W_{p+q} \rightarrow A_{m+n}$ is homogeneous of biweight $(p+q,m+n)$ (see Definition \ref{defcup}). Moreover, $b_K$ is homogeneous of biweight $(1,1)$. Thus the unital associative $k$-algebras $Hom(W_{\bullet},A)$ and $HK^{\bullet}(A)$ 
are $\mathbb{N}\times \mathbb{N}$-graded by the biweight. The homogeneous component of biweight $(p,m)$ of $HK^{\bullet}(A)$ is denoted $HK^p(A)_m$. Since
$$\partial_{\smile}: HK^p(A)_m \rightarrow HK^{p+1}(A)_{m+1},$$
the algebra $HK_{hi}^{\bullet}(A)$ is $\mathbb{N}\times \mathbb{N}$-graded by the biweight, and its $(p,m)$-component is denoted by $HK_{hi}^p(A)_m$. From Proposition \ref{zerohikoco}, we deduce the following.
\Bpo  \label{zerohikocoA}
Let $A=T(V)/(R)$ be a quadratic algebra. Assume that $V$ is finite dimensional. If $V\neq 0$, then $HK_{hi}^0(A)_0=0$. If $A$ is finite dimensional, $HK_{hi}^0(A)_{\max}=A_{\max}$ 
where $\max$ is the highest nonnegative integer $m$ such that $A_m\neq 0$. If the algebra $A$ is commutative, then for any $m\geq 0$, $HK_{hi}^0(A)_m$ equals the space of elements of $A_m$ annihilated by $V$.
\Epo

\subsection{Higher Koszul cohomology of symmetric algebras} \label{subsechosymalg}

Throughout this subsection, $A=S(V)$ is the symmetric algebra of the $k$-vector space $V$. We need no assumption on $\dim(V)$ or $\mathrm{char}(k)$. The following is standard.
\Blm
Let $V$ be a $k$-vector space and $A=S(V)$ be the symmetric algebra of $V$. For any $p\geq 0$, the space $W_p$ is equal to the image of the $k$-linear map
$Ant:V^{\otimes p}\rightarrow V^{\otimes p}$ defined by
$$Ant(v_1, \ldots , v_p)= \sum _{\sigma \in \Sigma_p} \mathrm{sgn}(\sigma)\, v_{\sigma (1)} \ldots v_{\sigma (p)},$$
for any $v_1, \ldots, v_p$ in $V$, where $\Sigma_p$ is the symmetric group and $\mathrm{sgn}$ is the signature.
\Elm
\Bpo \label{symalg}
Let $V$ be a $k$-vector space and $A=S(V)$ be the symmetric algebra of $V$. Let $M$ be a symmetric $A$-bimodule. 
The differentials $b_K$ of the complexes $M\otimes W_{\bullet}$ and $Hom(W_{\bullet},M)$ vanish. Therefore, $HK_{\bullet}(A,M)=M\otimes W_{\bullet}$ and $HK^{\bullet}(A,M)=Hom(W_{\bullet},M)$.
\Epo
\Bdm
Equation (\ref{defb}) can be written as
$$b_K(m \otimes x_1 \ldots x_{p}) =m.(x_1\otimes x_2 \ldots x_{p} +(-1)^p x_p  \otimes x_1 \ldots x_{p-1}),$$
and the right hand side vanishes according to the previous lemma and the relation
$$Ant(v_p, v_1, \ldots , v_{p-1})= (-1)^{p-1} Ant(v_1, \ldots , v_p).$$
Similarly, $b_K(f)=0$ for any Koszul cochain $f$.
\Edm
\\

Let us recall some facts about quadratic algebras~\cite{popo:quad}. Applying the functor $-\otimes_A k$ to the bimodule complex $K(A)=(A\otimes W_{\bullet} \otimes A, d)$, one obtains the left Koszul complex $K_{\ell}(A)=(A\otimes W_{\bullet}, d_{\ell})$ of left $A$-modules. The algebra $A$ is Koszul if and only if $K_{\ell}(A)$ is a 
resolution of $k$. Note that $\mu \otimes_A k$ coincides with the augmentation $\epsilon$. From (\ref{defd}) and using obvious notation, we have
\begin{equation} \label{defleftd}
d_{\ell}(a \otimes x_1 \ldots x_{p}) =ax_1\otimes x_2 \ldots x_{p}.
\end{equation}
\Bte \label{hkcosymalg}
Let $V$ be a $k$-vector space and $A=S(V)$ be the symmetric algebra of $V$. Assume that $\dim(V)=n$ is finite. We have
\begin{eqnarray*}
  HK_{hi}^n(A) &\cong& k,\\
  HK_{hi}^p(A) &\cong& 0 \ \mathrm{if} \ p\neq n.
\end{eqnarray*}
\Ete
\Bdm
Proposition \ref{symalg} shows that the differential $\partial_{\smile}$ on $HK^{\bullet}(A)$ coincides with the differential $e_A\underset{K}{\smile} -$ on $Hom(W_{\bullet},A)$. 
Given $f:W_p \rightarrow A$, denote by $F:A\otimes W_p\rightarrow A$ the left $A$-linear extension of $f$ to $A\otimes W_p$. From equation (\ref{defleftd}) applied to 
$1\otimes x_1 \ldots x_{p+1}$, and from
$$(e_A\underset{K}{\smile} f) (x_1 \ldots x_{p+1}) = (-1)^p x_1 . f(x_2 \ldots  x_{p+1}),$$
we deduce that
$$e_A\underset{K}{\smile} f= (-1)^p F \circ d_{\ell},$$
where $ d_{\ell}: A\otimes W_{p+1} \rightarrow A\otimes W_p$ is restricted to $W_{p+1}$. Thus the differential $e_A\underset{K}{\smile} -$ coincides with the opposite of the differential $Hom_A(d_{\ell}, A)$. Since $A$ is Koszul, we have obtained that
$$HK_{hi}^{\bullet}(A)\cong Ext^{\bullet}_A(k,A).$$
Using that $A$ is AS-Gorenstein of global dimension $n$~\cite{popo:quad}, the theorem is proved.
\Edm

\setcounter{equation}{0}

\section{The Koszul cap products}

\subsection{Definition and first properties}

As for the cup product, we define $\underset{K}{\frown}$ by restricting the usual $\frown$ and using the notation of Subsection 2.1.
\Bdf \label{defcap}
Let $A=T(V)/(R)$ be a quadratic algebra. Let $M$ and $P$ be $A$-bimodules. For any Koszul $p$-cochain $f:W_p\rightarrow P$ and any Koszul $q$-chain $z=m \otimes x_1 \ldots x_{q}$ 
in $M\otimes W_q$, we define the Koszul $(q-p)$-chains $f \underset{K}{\frown} z$ and $z \underset{K}{\frown} f$ with coefficients in $P\otimes_A M$ and $M\otimes_A P$ respectively, 
by the following equalities
\begin{equation} \label{kolcap}
f\underset{K}{\frown} z = (-1)^{(q-p)p} (f(x_{q-p+1} \ldots x_q)\otimes_A m) \otimes x_1 \ldots  x_{q-p},
\end{equation}
\begin{equation} \label{korcap}
  z\underset{K}{\frown} f = (-1)^{pq} (m\otimes_A f(x_1 \ldots x_p)) \otimes x_{p+1} \ldots  x_{q}.
\end{equation}
The element $f \underset{K}{\frown} z$ is called the left Koszul cap product of $f$ and $z$, while $z \underset{K}{\frown} f$ is called their right Koszul cap product. \Edf

If $q<p$, then one has $f \underset{K}{\frown} z=z \underset{K}{\frown} f=0$. By definition, we have
 \begin{equation} \label{hokolcap}
\tilde{\chi} (\chi^{\ast}(F) \underset{K}{\frown} z) = F \frown \tilde{\chi}(z),
\end{equation}
\begin{equation} \label{hokorcap}
\tilde{\chi} (z \underset{K}{\frown} \chi^{\ast}(F)) = \tilde{\chi}(z) \frown F,
\end{equation}
for any  Hochschild cochain $F:A^{\otimes p} \rightarrow P$ and any Koszul chain $z \in M\otimes W_q$. Considering both Koszul cap products $\underset{K}{\frown}$ respectively as left or right action, $M\otimes W_{\bullet}$ becomes a graded bimodule over the graded algebra $(Hom(W_{\bullet},A),\underset{K}{\smile})$, since these properties hold for the usual cup and cap products.

Similarly, we deduce the identities
\begin{equation} \label{koldga}
b_K(f \underset{K}{\frown} z)= b_K(f) \underset{K}{\frown} z +(-1)^p f \underset{K}{\frown} b_K(z),
\end{equation}
\begin{equation} \label{kordga}
b_K(z \underset{K}{\frown} f)= b_K(z) \underset{K}{\frown} f +(-1)^q z \underset{K}{\frown} b_K(f),
\end{equation}
from the identities known for the usual $\frown$. So 
$M\otimes W_{\bullet}$ is a differential graded bimodule over the dga $Hom(W_{\bullet},A)$. The proof of the following is clear.
\Bpo
Let $A=T(V)/(R)$ be a quadratic algebra. Both Koszul cap products $\underset{K}{\frown}$ at the chain-cochain level define Koszul cap products, still denoted by $\underset{K}{\frown}$, 
on Koszul (co)homology classes. Formulas (\ref{hokolcap}) and (\ref{hokorcap}) pass to classes. Considering Koszul cap products as actions, for any $A$-bimodule $M$, $HK_{\bullet}(A,M)$ is a graded bimodule on the graded algebra $HK^{\bullet}(A)$. 
In particular, $HK_{\bullet}(A,M)$ is a $Z(A)$-bimodule. Moreover, $HK_{\bullet}(A,k)=W_{\bullet}$ is a graded bimodule on the graded algebra $HK^{\bullet}(A,k)=W^{\ast}_{\bullet}$.
\Epo

\subsection{The Koszul cap bracket}

\Bdf \label{defcapbra}
Let $A=T(V)/(R)$ be a quadratic algebra. Let $M$ and $P$ be $A$-bimodules such that $M$ or $P$ is equal to $A$. For any Koszul $p$-cochain $f:W_p\rightarrow P$ and any Koszul $q$-chain 
$z \in M\otimes W_q$, we define the Koszul cap bracket $[f, z]_{\underset{K}{\frown}}$ by
\begin{equation} \label{kocapbra}
[f, z]_{\underset{K}{\frown}} =f\underset{K}{\frown} z - (-1)^{pq} z\underset{K}{\frown} f.
\end{equation}
\Edf

For $z=m \otimes x_1 \ldots x_{q}$, the explicit expression of the bracket is
\begin{equation} \label{kocapbraexp}
[f, z]_{\underset{K}{\frown}} = (-1)^{(q-p)p} f(x_{q-p+1} \ldots x_q)m \otimes x_1 \ldots  x_{q-p} - mf(x_1 \ldots x_p) \otimes x_{p+1} \ldots  x_{q}.
\end{equation}
If $p=0$, then $[f, z]_{\underset{K}{\frown}}= [f(1),m]_c\otimes x_1 \ldots x_{q}$, where $[-,-]_c$ denotes 
the commutator. The Koszul cap bracket passes to (co)homology classes. We still use the notation $[\alpha, \gamma]_{\underset{K}{\frown}}$ for classes 
$\alpha$ and $\gamma$ corresponding to  $f$ and $z$. When $M=A$, the maps $[f, -]_{\underset{K}{\frown}}$ and $[\alpha, -]_{\underset{K}{\frown}}$ are graded derivations of the 
graded $Hom(W_{\bullet},A)$-bimodule $A\otimes W_{\bullet}$, and of the graded $HK^{\bullet}(A)$-bimodule $HK_{\bullet}(A)$, respectively.

Similarly to what happens in cohomology, \emph{the Koszul differential $b_K$ in homology may be defined from the Koszul cap product}, and defining $b_K$ by (\ref{fundaho}) below, 
we may deduce the identities (\ref{koldga}) and (\ref{kordga}) from the derivation $[f, -]_{\underset{K}{\frown}}$. 
The subsequent theorem is analogous to Theorem \ref{thfundacoho}. The proof is left to the reader.
\Bte \label{thfundaho}
Let $A=T(V)/(R)$ be a quadratic algebra. For any Koszul cochain $z$ with coefficients in any $A$-bimodule $M$, we have the formula
\begin{equation} \label{fundaho}
[e_A, z]_{\underset{K}{\frown}}=-b_K(z).
\end{equation}
\Ete

\subsection{Actions of Koszul derivations} \label{subsecakd}

Using Subsection \ref{subseckoder}, we associate to a bimodule $M$ and a Koszul derivation $f: V \rightarrow M$ the derivation $D_f: A \rightarrow M$. 
The linear map 
$D_f\otimes Id_{W_{\bullet}}$ from $A\otimes W_{\bullet}$ to $M\otimes W_{\bullet}$
will still be denoted by $D_f$. The proof of the following proposition is easy.
\Bpo
Let $A=T(V)/(R)$ be a quadratic algebra and let $M$ be an $A$-bimodule. For any Koszul derivation $f: V \rightarrow M$ and any Koszul $q$-cycle $z \in A\otimes W_q$, 
\begin{equation} \label{kocapbraspe}
[f, z]_{\underset{K}{\frown}} =b_K(D_f(z)).
\end{equation}
\Epo
\Bcr \label{kocapcommutation}
Let $A=T(V)/(R)$ be a quadratic algebra and let $M$ be an $A$-bimodule. For any $p \in \{0,1,q\}$, $\alpha \in HK^p(A,M)$ and  $\gamma \in HK_q(A)$,
\begin{equation} \label{kocapbrazero}
[\alpha, \gamma]_{\underset{K}{\frown}} =0.
\end{equation}
\Ecr
\Bdm
The case $p=1$ follows from the proposition. The case $p=0$ is clear. Assume that $p=q$, $\alpha$ is the class of $f$ and $\gamma$ is the class of $z=a\otimes x_1\ldots x_p$. 
The equality (\ref{kocapbraexp}) gives
$$[f, z]_{\underset{K}{\frown}} = f(x_1 \ldots x_p).a - a.f(x_1 \ldots x_p)$$
which is an element of $[M,A]_c$. 
Since $[\alpha, \gamma]_{\underset{K}{\frown}}$ belongs to $HK_0(A,M)$, we conclude from the isomorphism
$$H(\tilde{\chi})_0 :HK_0(A,M)\rightarrow HH_0(A,M)=M/[M,A]_c.\ \Edm$$

Note that the same proof shows that $[\alpha, \gamma]_{\underset{K}{\frown}}=0$ if $\alpha \in HK^p(A)$ and $\gamma \in HK_p(A,M)$. We do not know whether the 
identity $[\alpha, \gamma]_{\underset{K}{\frown}}=0$ in the previous corollary holds for any $p$ and $q$ -- even if $A$ is Koszul! It holds for $M=A$ in the non-Koszul example of Section 9.

\setcounter{equation}{0}

\section{Higher Koszul homology} \label{sechkoho}

\subsection{Higher Koszul homology associated to a Koszul derivation} \label{subsecgenhkh}

A similar procedure to the one developed in Subsection \ref{subsechkoco} leads to the definition of a higher homology theory in the following situation. 
Let $A=T(V)/(R)$ be a quadratic algebra, $f: V \rightarrow A$ a Koszul derivation of $A$ and $M$ an $A$-bimodule. 
Assuming $\mathrm{char}(k)\neq 2$, the identity $[f]\underset{K}{\smile}[f]=0$ shows that the linear map $[f]\underset{K}{\frown} -$ is a chain 
differential on $HK_{\bullet}(A,M)$. 
We obtain therefore a new homology, called \emph{higher Koszul homology associated to $f$}. Analogously, $[D_f ] \frown -$ is a chain differential 
on $HH_{\bullet}(A,M)$, hence a \emph{higher Hochschild homology associated to $f$}. 
The map $H(\tilde{\chi})$ induces a morphism from the higher Koszul homology to the higher Hochschild homology, which is an isomorphism whenever $A$ is Koszul. For $z=m\otimes a_1\ldots a_p$ in $M\otimes A^{\otimes p}$, we deduce from the Hochschild analogue of equality (\ref{kolcap}) that 
$$D_f \frown z = (-1)^{p-1} (D_f(a_p) m) \otimes a_1 \ldots  a_{p-1}.$$
Thus \emph{$D_f \frown -$ coincides with the Rinehart-Goodwillie operator associated to the derivation $D_f$ of $A$}~\cite{rine:difform, good:cychom}.

\subsection{Higher Koszul homology associated to $e_A$} \label{subsechkh}

Let us fix $f=e_A$ for the rest of the paper. Without any assumption on the characteristic of $k$, the $k$-linear map $e_A\underset{K}{\frown} -$ is a chain differential 
on $M \otimes W_{\bullet}$, and next $\overline{e}_A\underset{K}{\frown} -$ is a chain differential on $HK_{\bullet}(A,M)$.

\Bdf \label{hikoho}
Let $A=T(V)/(R)$ be a quadratic algebra and let $M$ be an $A$-bimodule. The differential $\overline{e}_A\underset{K}{\frown} -$ of $HK_{\bullet}(A,M)$ will be denoted by 
$\partial_{\frown}$. The homology of $HK_{\bullet}(A,M)$ endowed with $\partial_{\frown}$ is called the higher Koszul homology of $A$ with coefficients in $M$ and is denoted 
by $HK^{hi}_{\bullet}(A,M)$. We set $HK^{hi}_{\bullet}(A)=HK^{hi}_{\bullet}(A,A)$.
\Edf

If we want to evaluate $\partial_{\frown}$ on classes, it suffices to go back to the formula
$$e_A\underset{K}{\frown} z = mx_1 \otimes x_2 \ldots  x_p$$
for any cycle $z=m\otimes x_1\ldots x_p$ in $M \otimes W_p$. If $M=k$, 
the differential $e_A\underset{K}{\frown} -$ vanishes, so $HK^{hi}_p(A,k)=W_p$ for any $p\geq 0$.

\subsection{Higher Koszul homology with coefficients in $A$} \label{subsechkhAcoeff}

\Blm
Let $A=T(V)/(R)$ be a quadratic algebra. Given $\alpha$ in $HK^p(A)$ and $\gamma$ in $HK_q(A)$, the following equalities hold
$$\partial_{\frown}(\alpha \underset{K}{\frown} \gamma)= \partial_{\smile}(\alpha) \underset{K}{\frown} \gamma = (-1)^p \alpha \underset{K}{\frown} \partial_{\frown}(\gamma),$$
$$\partial_{\frown}(\gamma \underset{K}{\frown} \alpha)= \partial_{\frown}(\gamma) \underset{K}{\frown} \alpha = (-1)^q \gamma \underset{K}{\frown} \partial_{\smile}(\alpha).$$
\Elm

The proof is left to the reader. Consequently, the Koszul cap products are defined in $HK_{hi}^{\bullet}(A)$ acting on $HK^{hi}_{\bullet}(A)$ and are still denoted by $\underset{K}{\frown}$. 
This makes $HK^{hi}_{\bullet}(A)$ a graded bimodule over the graded algebra $HK_{hi}^{\bullet}(A)$. More generally, $HK^{hi}_{\bullet}(A,M)$ is a graded bimodule over the graded 
algebra $HK_{hi}^{\bullet}(A)$ for any $A$-bimodule $M$.

As we have already done in cohomology, but without any assumption on $V$, we show that the space $HK^{hi}_{\bullet}(A)$ is bigraded. A Koszul $q$-chain $z$ in $A_n \otimes W_q$ is said to be homogeneous of weight $n$. The space of Koszul 
chains $A \otimes W_{\bullet}$ is $\mathbb{N}\times \mathbb{N}$-graded by the \emph{biweight} $(q,n)$, where $q$ is called the \emph{homological weight} and $n$ is called the 
\emph{coefficient weight}. Moreover, $b_K$ is homogeneous of biweight $(-1,1)$. Thus the space $HK_{\bullet}(A)$ is $\mathbb{N}\times \mathbb{N}$-graded by the biweight. 
The homogeneous component of biweight $(q,n)$ of $HK_{\bullet}(A)$ is denoted by $HK_q(A)_n$. Since 
$$\partial_{\frown}: HK_q(A)_n \rightarrow HK_{q-1}(A)_{n+1},$$
the space $HK^{hi}_{\bullet}(A)$ is $\mathbb{N}\times \mathbb{N}$-graded by the biweight, and its $(q,n)$-component is denoted by $HK^{hi}_q(A)_n$.

Assume now that $V$ is finite dimensional. If $f:W_p \rightarrow A_m$ and $z \in A_n \otimes W_q$ are homogeneous of biweights $(p,m)$ and $(q,n)$ respectively, 
then $f\underset{K}{\frown} z$ and $z\underset{K}{\frown} f$ are homogeneous of biweight $(q-p, m+n)$ where
\begin{equation} \label{kolcapA}
f\underset{K}{\frown} z = (-1)^{(q-p)p} f(x_{q-p+1} \ldots x_q)a \otimes x_1 \ldots  x_{q-p},
\end{equation}
\begin{equation} \label{korcapA}
  z\underset{K}{\frown} f = (-1)^{pq} a f(x_1 \ldots x_p) \otimes x_{p+1} \ldots  x_{q},
\end{equation}
and $z=a\otimes x_1\ldots x_q$. The $Hom(W_{\bullet},A)$-bimodule $A \otimes W_{\bullet}$, the $HK^{\bullet}(A)$-bimodule $HK_{\bullet}(A)$ and the 
$HK_{hi}^{\bullet}(A)$-bimodule $HK^{hi}_{\bullet}(A)$ are thus $\mathbb{N}\times \mathbb{N}$-graded by the biweight. The proof of the following is left to the reader.
\Bpo  \label{easyhikohoA}
For any quadratic algebra $A=T(V)/(R)$,
$$HK_0(A)_0=HK^{hi}_0(A)_0=k.$$
Moreover $HK_0(A)_1= HK_1(A)_0=V$ and $\partial_{\frown}: HK_1(A)_0 \rightarrow HK_0(A)_1$ is the identity map of $V$. 
As a consequence,
$$HK^{hi}_0(A)_1=HK^{hi}_1(A)_0=0.$$
\Epo

\subsection{Higher Koszul homology of symmetric algebras}

\Bte \label{hkhsymalg}
Given a $k$-vector space $V$ and the symmetric algebra $A=S(V)$, we have 
\begin{eqnarray*}
  HK^{hi}_0(A) &\cong& k,\\
  HK^{hi}_p(A) &\cong& 0 \ \mathrm{if} \ p>0.
\end{eqnarray*}
\Ete
\Bdm
Proposition \ref{symalg} shows that the differential $\partial_{\frown}$ on $HK_{\bullet}(A)$ coincides with the differential $e_A\underset{K}{\frown} -$ on $A\otimes W_{\bullet}$. 
From equation (\ref{defleftd}), we see that the complex $(HK_{\bullet}(A),\partial_{\frown})$ coincides with the left Koszul complex $K_{\ell}(A)=(A\otimes W_{\bullet}, d_{\ell})$. Since $A$ is Koszul, we deduce $HK^{hi}_{\bullet}(A)$ as stated. 
\Edm
\\

Our aim is now to generalize this theorem to any Koszul algebra, in characteristic zero. This generalization is presented in the next section. The proof given below uses some 
standard facts on Hochschild homology of graded algebras including the Rinehart-Goodwillie operator. 

\setcounter{equation}{0}

\section{Higher Koszul homology and de Rham cohomology} \label{secderham}

\subsection{Standard facts on Hochschild homology of graded algebras} \label{subsecstandard}

For Hochschild homology of graded algebras, we refer to Goodwillie~\cite{good:cychom}, Section 4.1 of Loday's book~\cite{loday:cychom}, or Section 9.9 of Weibel's book~\cite{weib:homo}. 
In this subsection, $A$ is a unital associative $k$-algebra which is $\mathbb{N}$-graded by a weight. The homogeneous component of weight $p$ of $A$ is denoted by $A_p$ and 
we set $|a|=p$ for any $a$ in $A_p$. We assume that $A$ is connected, i.e. $A_0=k$, so that $A$ is augmented. Recall that the weight map $D=D_A: A \rightarrow A$ of the graded algebra $A$ is defined by
$D(a)=pa$
for any $p\geq 0$ and $a$ in $A_p$. As recalled in Subsection \ref{subsecgenhkh}, the Rinehart-Goodwillie operator $e_D=D \frown -$ of $A\otimes A^{\otimes \bullet}$ is defined by
 $$e_D(a\otimes a_1\ldots a_p) = (-1)^{p-1} (|a_p|a_pa) \otimes a_1 \ldots  a_{p-1},$$
for any $a$, $a_1, \ldots, a_p$ in $A$ with $a_p$ homogeneous. If $p=0$, note that $e_D(A)=0$.

Denote by $[D]$ the Hochschild cohomology class of $D$. Assuming $char(k)\neq 2$, Gerstenhaber's identity $2 D \smile D = b(D \circ D)$ shows that the map $H(e_D)=[D] \frown -$ is a chain differential on $HH_{\bullet}(A)$, and $[D] \smile -$ is a cochain differential on $HH^{\bullet}(A)$. 
We denote by $HH^{hi}_{\bullet}(A)$ (resp. $HH_{hi}^{\bullet}(A)$) the so-obtained higher Hochschild homology (resp. cohomology) of $A$ 
with coefficients in $A$, already defined if $A$ is a quadratic algebra in Subsections 3.5 and 5.1.

Let $B$ be the normalized Connes differential of $A\otimes \bar{A}^{\otimes \bullet}$ where $\bar{A}=A/k$~\cite{loday:cychom,weib:homo}. Denoting the augmentation of $A$ by $\epsilon $, we identify $\bar{A}$ to the subspace $\ker(\epsilon)=\bigoplus_{m>0} A_m$ of $A$. Recall that
\begin{equation} \label{connesB}
B(a\otimes a_1\ldots a_p)= \sum_{0\leq i \leq p} (-1)^{pi} 1\otimes (a_{p-i+1}\ldots a_p\bar{a}a_1\ldots a_{p-i}),
\end{equation}
for any $a \in A$, and $a_1, \ldots, a_p$ in $\bar{A}$, where $\bar{a}$ denotes the class of $a$ in $\bar{A}$. Note that $B(a)=1\otimes \bar{a}$ for any $a$ in $A$. 
The operator $B$ passes to Hochschild homology and defines the cochain differential $H(B)$ on $HH_{\bullet}(A)$. We follow Van den Bergh~\cite{vdb:hom} for the subsequent definition. 
\Bdf 
The complex $(HH_{\bullet}(A), H(B))$ is called the de Rham complex of $A$. The homology of this complex is called the de Rham cohomology of $A$ and is denoted by $H_{dR}^{\bullet}(A)$. 
\Edf

If $char(k)=0$, it turns out that \emph{one of both differentials $H(B)$ and $H(e_D)$ of $HH_{\bullet}(A)$ is -- up to a normalization -- a contracting 
homotopy of the other one}. This duality linking $H(B)$ and $H(e_D)$ is a consequence of the Rinehart-Goodwillie identity (\ref{rinegood}) below. Let us introduce the weight 
map $L_D$ of $A\otimes \bar{A}^{\otimes \bullet}$ by
$$L_D(z)=|z|z,$$
for any homogeneous $z=a\otimes a_1\ldots a_p$, where $|z|=|a|+ |a_1| + \cdots + |a_p|$. Clearly, $L_D$ defines an operator $H(L_D)$ on $HH_{\bullet}(A)$. 
Note that $A\otimes A^{\otimes \bullet}$, $HH_{\bullet}(A)$ and $HH^{hi}_{\bullet}(A)$ are graded by the \emph{total weight} (called simply the weight), and that the 
operators $H(e_D)$, $H(B)$ and $H(L_D)$ are weight homogeneous. Let us state the Rinehart-Goodwillie identity; for a proof, see for example Corollary 4.1.9 in~\cite{loday:cychom}.
\Bpo
Let $A$ be a connected $\mathbb{N}$-graded $k$-algebra. The identity  \begin{equation} \label{rinegood}
[H(e_D), H(B)]_{gc} =H(L_D),
\end{equation}
holds, where $[-,-]_{gc}$ denotes the graded commutator with respect to the homological degree.
\Epo

The following consequence is a noncommutative analogue of Poincar\'e Lemma.
\Bte \label{poincarelemma}
Let $A$ be a connected $\mathbb{N}$-graded $k$-algebra. Assume $char(k)=0$. We have
\begin{eqnarray*}
  H_{dR}^0(A)\ \cong \  HH^{hi}_0(A)  &\cong & k,\\
  H_{dR}^p(A)\ \cong \  HH^{hi}_p(A)  & \cong & 0 \ \mathrm{if} \ p>0.
\end{eqnarray*}
\Ete
\Bdm
Let $\alpha \neq 0$ be a weight homogeneous element in $HH_p(A)$. Assume that $H(e_D)(\alpha)=0$. The identity (\ref{rinegood}) provides
\begin{equation} \label{rinegood2}
H(e_D)\circ H(B)(\alpha) = |\alpha| \alpha.
\end{equation}
If $p>0$, then $|\alpha|\neq 0$, so that $\alpha$ is a $H(e_D)$-boundary, showing that $HH^{hi}_p(A)=0$. If $p=0$, any $\alpha$ in $HH_0(A)$ is a cycle for $H(e_D)$ and if 
$|\alpha|\neq 0$, it is a boundary by (\ref{rinegood2}). If $p=|\alpha|= 0$, $\alpha$ cannot be a boundary since $H(e_D)$ adds 1 to the coefficient weight. Thus $HH^{hi}_0(A)=k$. The proof for the de Rham case is similar. Note that the assumption $char(k)=0$ is essential in this proof, except for proving that $H_{dR}^0(A)\cong k$ and that $HH^{hi}_0(A)_0\cong k$.
\Edm

\subsection{Consequences for quadratic algebras}

If $A$ is quadratic, then $H(\tilde{\chi}): HK_p(A) \rightarrow HH_p(A)$ is always an isomorphism for $p=0$ and $p=1$; moreover, if $A$ is Koszul 
it is an isomorphism for any $p$. As a consequence, $H(\tilde{\chi})$ induces an isomorphism from $HK^{hi}_p(A)$ 
to $HH^{hi}_p(A)$ for  $p=0$, and for any $p$ if $A$ is Koszul. So, generalizing Theorem \ref{hkhsymalg} in characteristic zero, we obtain 
the following consequence of the previous theorem.
\Bte \label{hkhkoszul}
Let $A=T(V)/(R)$ be a quadratic algebra. Assume that $char(k)=0$. We have $HK^{hi}_0(A) \cong k$. If $A$ is Koszul, then for any $p>0$, 
$$HK^{hi}_p(A) \cong 0.$$
\Ete

It would be more satisfactory to find a proof within the Koszul calculus, possibly without any assumption on $char(k)$. 
We would also like to know if the converse of this theorem holds, namely, if the following conjecture is true.
\Bcn  \label{conhkoszul}
Let $A=T(V)/(R)$ be a quadratic algebra. The algebra $A$ is Koszul if and only if there are isomorphisms 
\begin{eqnarray*}
  HK^{hi}_0(A)  &\cong & k,\\
  HK^{hi}_p(A)  & \cong & 0 \ \mathrm{if} \ p>0.
\end{eqnarray*}
\Ecn

Let us comment on this conjecture. In the non-Koszul example of Section 9, we will find that $HK^{hi}_2(A)\neq 0$ -- agreeing the conjecture. 
Within the graded Hochschild calculus, this conjecture is meaningless, since \emph{any} graded algebra has a trivial higher Hochschild homology 
as stated in Theorem \ref{poincarelemma}. Consequently, the higher Koszul homology provides more information on quadratic algebras than the 
higher Hochschild homology. Moreover, if Conjecture \ref{conhkoszul} is true, then \emph{the Koszul algebras would be exactly the acyclic objects for the
higher Koszul homology}.

In Subsection \ref{subsechosymalg}, the left Koszul complex $K_{\ell}(A)=K(A) \otimes_A k$ associated to any quadratic algebra $A$ was recalled. 
Since $A$ is Koszul if and only if $K_{\ell}(A)$ is a resolution of $k$, Conjecture \ref{conhkoszul} is an immediate consequence of the following.
\Bcn  \label{lconhkoszul}
Let $A=T(V)/(R)$ be a quadratic algebra. For any $p\geq 0$ \begin{equation}
  HK^{hi}_p(A)  \cong H_p(K_{\ell}(A)).
\end{equation}
\Ecn

A stronger conjecture asserts that there exists a quasi-isomorphism from the complex $(HK_{\bullet}(A),\partial_{\frown})$ to the 
complex $K_{\ell}(A)$. The proof of Theorem \ref{hkhsymalg} shows that the stronger conjecture holds for symmetric algebras. 
For any quadratic algebra $A$, it is well-known that $H_0(K_{\ell}(A))\cong k$ and $H_1(K_{\ell}(A))\cong 0$, therefore Conjecture 
\ref{lconhkoszul} would imply that $HK^{hi}_0(A)\cong k$ and $HK^{hi}_1(A)\cong 0$. What we know about $HK^{hi}_1(A)$ is 
that $HK^{hi}_1(A)_0\cong 0$ (Proposition \ref{easyhikohoA}), and $HK^{hi}_1(A)_1\cong 0$ (next subsection). 
Note that the non-Koszul example of Section 9 will satisfy Conjecture \ref{lconhkoszul}.

\subsection{The Connes differential on Koszul classes}

From the equality (\ref{connesB}) defining the Connes differential $B$ of $A\otimes \bar{A}^{\otimes \bullet}$, observe that $B(A\otimes W_p)$ is 
not included in $A\otimes W_{p+1}$, so that it seems hard to find an analogue to $B$ at the Koszul chain level. We prefer to search 
an analogue to $H(B)$ at the Koszul homology level. In this subsection, the notation $H(B)$ is simplified and replaced by $B$. We are interested in 
the following question. 
Let $A=T(V)/(R)$ be a quadratic algebra. 

Does there exist a $k$-linear cochain differential $B_K$ on $HK_{\bullet}(A)$ such that the diagram
\begin{eqnarray} \label{question}
HK_p(A)  \stackrel{B_K}{\longrightarrow} & HK_{p+1}(A) \nonumber  \\
\downarrow H(\tilde{\chi})_p  \ \ \ \ \  &  \downarrow H(\tilde{\chi})_{p+1} \\
HH_p(A) \stackrel{B}{\longrightarrow} & HH_{p+1}(A) \nonumber \end{eqnarray}
commutes for any $p\geq 0$? 

Since $B$ and $H(\tilde{\chi})$ preserve the total weight, $B_K$ should preserve the total weight too. 
Therefore, using our notation for coefficient weight, we impose that
$$B_K : HK_p(A)_m \rightarrow HK_{p+1}(A)_{m-1}.$$

The answer to the question is affirmative if $A$ is Koszul since the vertical arrows are isomorphisms, and in this case the 
corresponding Rinehart-Goodwillie identity linking the differentials $B_K$ and $\partial_{\frown}$ of $HK_{\bullet}(A)$ holds. 
If the answer is affirmative for a non-Koszul algebra $A$, Conjecture \ref{conhkoszul} would imply that this Koszul Rinehart-Goodwillie identity 
does not hold in characteristic zero, and it would be interesting to measure the defect to be an identity, e.g. in the explicit example of Section 9.

Let us begin by examining the diagram (\ref{question}) for $p=0$. In this case, such a $B_K$ exists since the vertical arrows are isomorphisms.
It suffices to pre and post compose the map
$$B:HH_0(A)\rightarrow HH_1(A), \  [a] \mapsto [1\otimes \bar{a}],$$
with the isomorphism and its inverse in order to obtain $B_K$, however an explicit expression of $B_K$ is not clear. 
It is easy to obtain it for small coefficient weights. Clearly,
$$B_K: HK_0(A)_1=V \rightarrow HK_1(A)_0=V$$
is the identity of $V$. 
Next, assume $char(k)\neq 2$ and consider the projections $ant$ and $sym$ of $V\otimes V$ defined by
$$ant(x\otimes y)=\frac{1}{2}(x\otimes y - y\otimes x), \ sym(x\otimes y)=\frac{1}{2}(x\otimes y + y\otimes x),$$
for any $x$ and $y$ in $V$. 
The proof of the following lemma is straightforward.
\Blm
Let $A=T(V)/(R)$ be a quadratic algebra. If $char(k)\neq 2$, we have 
$$HK_2(A)_0 = R \cap ant(V\otimes V), \ HK_1(A)_1=\frac{ant^{-1}(R)}{sym(R)}, \ HK_0(A)_2 \cong \frac{V\otimes V}{ant(V\otimes V)+R}\,.$$
\Elm

The map $B_K: HK_0(A)_2 \rightarrow HK_1(A)_1$ is thus defined by $B_K([a])=[sym(a)]$ for any $[a]$ in $\frac{V\otimes V}{ant(V\otimes V)+R}$. Let us continue a bit further by defining the map
$$B_K: HK_1(A)_1 \rightarrow HK_2(A)_0$$
by $B_K([a])=2 ant(a)$ for any $[a]$ in $\frac{ant^{-1}(R)}{sym(R)}$. 
The proof of the following lemma is direct.
\Blm \label{tw2}
The map $B_K: HK_1(A)_1 \rightarrow HK_2(A)_0$ is surjective and together with $B_K: HK_0(A)_2 \rightarrow HK_1(A)_1$ it satisfies 
the Koszul Rinehart-Goodwillie identity
$$(\partial_{\frown} \circ B_K + B_K\circ \partial_{\frown})([a])=2[a],$$
for any $[a]$ in $HK_1(A)_1$. Moreover, $H(\tilde{\chi})_2:HK_2(A)_0 \rightarrow HH_2(A)_2$ is an isomorphism. 
\Elm

Note that $HH_p(A)_t$ denotes the homogeneous component of \emph{total weight} $t$. Using the previous $B_K$, the diagram (\ref{question}) corresponding to $p=1$ and total weight 2 commutes. From Lemma \ref{tw2}, we obtain immediately the following proposition.
\Bpo
Let $A=T(V)/(R)$ be a quadratic algebra. If $char(k)\neq 2$, we have
$$HK^{hi}_2(A)_0 \cong HK^{hi}_1(A)_1 \cong 0.$$
\Epo

Generalizing $B_K: HK_1(A)_1 \rightarrow HK_2(A)_0$ as below, we obtain the following.
\Bpo \label{p0}
Let $A=T(V)/(R)$ be a quadratic algebra. If $p\geq 2$ is not divisible by $char(k)$, then $HK^{hi}_p(A)_0 \cong 0$.
\Epo
\Bdm
Denote $b_{K,p}: W_p \rightarrow V \otimes W_{p-1}$ and $b_{K,p-1}: V \otimes W_{p-1} \rightarrow A_2 \otimes W_{p-2}$ the differential $b_K$ 
on $p$-chains of weight 0 and on $(p-1)$-chains of weight 1. We have
$$HK_p(A)_0= \ker(b_{K,p})\subseteq W_p \subseteq V\otimes  W_{p-1}, \ \  \ HK_{p-1}(A)_1=\frac{\ker (b_{K,p-1})}{\mathrm{im} (b_{K,p})},$$
and $e_A \frown z=z$ for any $z$ in $\ker(b_{K,p})$. The map
$$\partial_{\frown}: HK_p(A)_0 \rightarrow HK_{p-1}(A)_1$$
is defined by $\partial_{\frown}(z)=[z]$ for any $z$ in $\ker(b_{K,p})$. In order to show that this map is injective under the hypothesis on 
the characteristic, it suffices to define
$$B_K: HK_{p-1}(A)_1 \rightarrow HK_p(A)_0$$
such that $B_K \circ \partial_{\frown} = p\, \mathrm{Id}_{HK_p(A)_0}$.
For this, restrict the operators $t$ and $N$ of cyclic homology~\cite{loday:cychom} to $V^{\otimes p}$. We get the operators $\tau$ and $\gamma$ 
of $V^{\otimes p}$ given for any $v_1, \ldots, v_p$ in $V$ and $z$ in $V^{\otimes p}$ by
$$\tau (v_1 \otimes \ldots \otimes v_p)= (-1)^{p-1} v_p\otimes  v_1\otimes \ldots \otimes v_{p-1},$$
$$\gamma (z)= z+\tau(z) + \cdots + \tau ^{p-1}(z).$$
Clearly $\tau^p= \mathrm{Id}_{V^{\otimes p}}$ and $(1-\tau)\circ \gamma = \gamma \circ (1- \tau)=0$. 
\Blm
If $z\in V \otimes W_{p-1}$ is such that $b_{K,p-1}(z)=0$, then $\gamma (z) \in W_p$ and $b_{K,p}(\gamma (z))=0$.
\Elm
\Bdm
Write $z=x\otimes x_1\ldots x_{p-1}$ with usual notation. For $1 \leq i \leq p-1$, define
$$\mu_{i,i+1}= \mathrm{Id}_{V^{\otimes i-1}}\otimes \mu \otimes \mathrm{Id}_{V^{\otimes p-i-1}}:{V^{\otimes p}} \rightarrow {V^{\otimes i-1}}\otimes A_2 \otimes {V^{\otimes p-i-1}},$$
so that $\mu_{i,i+1}(v_1\otimes \ldots \otimes v_p)= v_1\otimes \ldots \otimes v_{i-1} \otimes (v_iv_{i+1}) \otimes \ldots \otimes v_p$. 
Clearly,
\begin{equation} \label{mutau}
\mu_{i+1,i+2}\circ \tau = -\tau \circ \mu_{i,i+1}
\end{equation}
where $\tau$ on the right-hand side acts on $A^{\otimes p-1}$ by the same formula, hence with sign $(-1)^{p-2}$. The formula
$$b_{K,p-1}(z)=(xx_1)\otimes x_2 \ldots x_{p-1} +(-1)^{p-1} (x_{p-1}x) \otimes x_1\ldots x_{p-2}$$
shows that $b_{K,p-1}$ coincides with the restriction of $\mu_{1,2}\circ (1+\tau)$ to $V \otimes W_{p-1}$. Since $\gamma (z)$ is equal to 
$$x \otimes x_1\ldots x_{p-1}+(-1)^{p-1}x_{p-1}\otimes x \otimes x_1 \ldots x_{p-2}+ x_{p-2}\otimes x_{p-1} \ldots x_{p-3} + \cdots 
+ (-1)^{p-1}x_1\otimes x_2\ldots x$$
we see that
$$\mu_{1,2}(\gamma(z))=\mu_{1,2}(z+\tau(z))=b_{K,p-1}(z)=0$$
by assumption. Therefore, using equation (\ref{mutau}), we get
$$\mu_{2,3}(\gamma(z))=\mu_{2,3}(\tau(z)+\tau^2(z))=-\tau \circ \mu_{1,2}(z+ \tau(z))=0,$$
and we proceed inductively, up to
$$\mu_{p-1,p}(\gamma(z))=\mu_{p-1,p}(\tau^{p-2}(z)+\tau^{p-1}(z))=-\tau \circ \mu_{p-2,p-1}(\tau^{p-3}(z)+ \tau^{p-2}(z))=0.$$
Thus, we have proved successively that $\gamma (z)$ belongs to $R\otimes V^{\otimes p-2}$, $V\otimes R\otimes V^{\otimes p-3}$, up 
to $V^{\otimes p-2} \otimes R$, which means that $\gamma (z) \in W_p$. Next, the equality $b_{K,p}(\gamma (z))=0$ is clear since $b_{K,p}$ 
coincides with the restriction of $1-\tau$ to $W_p$.
\Edm

\medskip

So we set $B_K([z])=\gamma (z)$ for any $[z]$ in $HK_{p-1}(A)_1$ where $z\in \ker (b_{K,p-1})$. It is immediate that 
$(B_K \circ \partial_{\frown})(z) = \gamma (z)=p z$ for any $z$ in $\ker(b_{K,p})$. Proposition \ref{p0} is thus proved. \Edm
\\

Note that the corresponding diagram (\ref{question}) w.r.t. $p-1$ and total weight $p$ commutes. 
Remark as well that $H_p(K_{\ell}(A))_0=0$, thus Conjecture \ref{lconhkoszul} is satisfied in characteristic zero for coefficient weight zero. 

\setcounter{equation}{0}

\section{Higher Koszul cohomology and Calabi-Yau algebras} \label{seccy}

For the definition of Calabi-Yau algebras, we refer to Ginzburg~\cite{vg:cy}. The following is a higher Hochschild cohomology version of 
Poincar\'e duality, and it is based on the material recalled in Subsection \ref{subsecstandard}.
\Bte \label{poincareduality}
Let $A$ be a connected $\mathbb{N}$-graded $k$-algebra. Assume that $char(k)=0$. If $A$ is $n$-Calabi-Yau, 
then  \begin{eqnarray*}
  HH_{hi}^n(A)  &\cong & k, \\
  HH_{hi}^p(A)  & \cong & 0 \ \mathrm{if} \ p\neq n.
\end{eqnarray*}
\Ete
\Bdm
Let $c\in HH_n(A)$ be the fundamental class of the Calabi-Yau algebra $A$. As proved by the second author in~\cite{tl:bvcy} (Th\'eor\`eme 4.2), 
the Van den Bergh duality~\cite{vdb:dual} can be expressed by saying that the $k$-linear map
$$-\frown c: HH^p(A,M) \longrightarrow HH_{n-p}(A,M)$$
is an isomorphism for any $p$ and any $A$-bimodule $M$. As in Subsection \ref{subsecstandard}, $D$ denotes the weight map of $A$, the map 
$[D] \frown -$ is a chain differential on $HH_{\bullet}(A)$, and $[D] \smile -$ is a cochain differential on $HH^{\bullet}(A)$. Clearly the diagram
\begin{eqnarray} \label{poincareduality}
HH^p(A)  \ \ \ \stackrel{[D] \smile -}{\longrightarrow} & HH^{p+1}(A) \nonumber  \\
\downarrow  -\frown c \ \ \ \ \ \ \ \  &  \downarrow -\frown c \\
HH_{n-p}(A) \ \ \ \stackrel{[D] \frown -}{\longrightarrow} & HH_{n-p-1}(A) \nonumber \end{eqnarray}
commutes for any $p\geq 0$. 
Since the vertical arrows are isomorphisms, they induce isomorphisms $HH_{hi}^p(A) \cong H^{hi}_{n-p}(A)$. 
The result thus follows from Theorem \ref{poincarelemma}.
\Edm

\Bcr \label{hkcohkoszul}
Let $A=T(V)/(R)$ be a quadratic algebra. Assume that $char(k)=0$. If $A$ is Koszul and $n$-Calabi-Yau, then
\begin{eqnarray*}
  HK_{hi}^n(A)  &\cong & k, \\
  HK_{hi}^p(A)  & \cong & 0 \ \mathrm{if} \ p\neq n.
\end{eqnarray*}
\Ecr
\Bdm
Since $A$ is Koszul, $H(\chi^{\ast})$ induces an isomorphism from $HH_{hi}^{\bullet}(A)$ to $HK_{hi}^{\bullet}(A)$.
\Edm
\\

Analogously to Conjecture \ref{conhkoszul}, we formulate the following.
\Bcn  \label{concohkoszul}
Let $A=T(V)/(R)$ be a Koszul quadratic algebra. The algebra $A$ is $n$-Calabi-Yau if and only if there are isomorphisms 
\begin{eqnarray*}
  HK_{hi}^n(A)  &\cong & k, \\
  HK_{hi}^p(A)  & \cong & 0 \ \mathrm{if} \ p\neq n.
\end{eqnarray*}
\Ecn

We will illustrate this conjecture by the example $A=T(V)$ when $\dim(V) \geq 2$. The complex $K_{\ell}(A)$ is in this case
$$0 \longrightarrow A\otimes V \stackrel{\mu}{\longrightarrow} A \longrightarrow 0,$$
so that $A$ is Koszul of global dimension $1$, and $A$ is not AS-Gorenstein since $\dim(V) \geq 2$, thus $A$ is not Calabi-Yau. 
The following proposition shows that Conjecture \ref{concohkoszul} is valid for these algebras.
\Bpo \label{tensoralgebra}
Let $V$ be a finite-dimensional $k$-vector space such that $\dim(V)\geq 2$, and $A=T(V)$ the tensor algebra of $V$. We have 
\begin{equation} \label{hkcohtensor}
\begin{array}{lll}
HK^0_{hi}(A) & \cong & 0   \\
HK^1_{hi}(A)_0 & \cong & V^{\ast}  \\
HK^1_{hi}(A)_1 & \cong & Hom(V, V)/k .\mathrm{Id}_V  \\
HK^1_{hi}(A)_m & \cong & Hom(V, V^{\otimes m})/<v\mapsto av-va; a\in V^{\otimes m-1})>  \ \ \mathrm{if} \ m\geq 2  \\
HK^p_{hi}(A) & \cong & 0 \ \mathrm{if} \ p\geq 2.  \end{array}
\end{equation}
\Epo
\Bdm
The homology of the complex $0 \longrightarrow A \stackrel{b_K}{\longrightarrow} Hom(V, A)  \longrightarrow 0$, where $b_K(a)(v)=av-va$ for any $a$ in $A$ and $v$ in $V$, is $HK^{\bullet}(A)$. Thus \begin{equation} \label{kcohtensor}
\begin{array}{lll}
HK^0(A) & \cong & Z(A) \ \cong \ k   \\
HK^1(A) & \cong & Hom(V, A)/<v\mapsto av-va; a\in A>  \\
HK^p(A) & \cong & 0 \ \mathrm{if} \ p\geq 2.  \end{array}
\end{equation}
Next, $\partial_{\smile}$ is defined from $HK^0(A)_0 \cong k$ to $HK^1(A)_1 \cong Hom(V,V)$ by $\partial_{\smile}(\lambda)=\lambda .\mathrm{Id}_V$
for any $\lambda$ in $k$, hence it is injective. Equations (\ref{hkcohtensor}) follow immediately. 
\Edm

\setcounter{equation}{0}

\section{Application of Koszul calculus to Koszul duality} \label{secduality}

Throughout this section, $V$ denotes a finite dimensional $k$-vector space and $A=T(V)/(R)$ is a quadratic algebra. Let $V^{\ast}= Hom(V,k)$ be 
the dual vector space of $V$. For any $p\geq 0$, the natural isomorphism from $(V^{\otimes p})^{\ast}$ to $V^{\ast \otimes p}$ is always 
understood \emph{without sign}. The reason is that in this paper, we are only interested in the \emph{ungraded situation}, meaning that there is 
no additional $\mathbb{Z}$-grading on $V$. Let $R^{\perp}$ be the subspace of $V^{\ast}\otimes V^{\ast}$ defined as the orthogonal of the subspace 
$R$ of $V\otimes V$, w.r.t. the natural duality between the space $V\otimes V$ and its dual $(V\otimes V)^{\ast}\cong V^{\ast}\otimes V^{\ast}$.
\Bdf \label{defdual}
The quadratic algebra $A^!=T(V^{\ast})/(R^{\perp})$ is called the Koszul dual of the quadratic algebra $A$.
\Edf
Recall that $A$ is Koszul if and only if $A^!$ is Koszul~\cite{popo:quad}. The homogeneous component of weight $m$ of $A^!$ is denoted by $A^!_m$.
The subspace of $V^{\ast \otimes p}$ corresponding to the subspace $W_p$ of $V^{\otimes p}$ is denoted by $W^!_p$. By definition, 
\begin{equation} \label{A!m}
A^!_m=V^{\ast \otimes m} / \sum_{i+2+j=m} V^{\ast \otimes i}\otimes R^{\perp} \otimes V^{\ast \otimes j},
\end{equation}
\begin{equation} \label{W!p}
W^!_p=\bigcap_{i+2+j=p} V^{\ast \otimes i}\otimes R^{\perp} \otimes V^{\ast \otimes j}.
\end{equation}

\subsection{Koszul duality in cohomology}  \label{subseckdco}

Recall that $HK^{\bullet}(A)$ is $\mathbb{N}\times \mathbb{N}$-graded by the biweight $(p,m)$, where $p$ is the homological weight and $m$ is 
the coefficient weight. The homogeneous component of biweight $(p,m)$ of $HK^{\bullet}(A)$ is denoted 
by $HK^p(A)_m$. It will be crucial for the Koszul duality to exchange the weights $p$ and $m$ in the definition of the Koszul cohomology of $A$, leading to a modified version of the Koszul cohomology algebra denoted by tilde accents. More precisely, for Koszul cochains 
$f:W_p \rightarrow A_m$ and $g:W_q \rightarrow A_n$, define $\tilde{b}_K(f)$ and $f \tilde{\underset{K}{\smile}} g$ by
\begin{equation} \label{tildedefcob}
\tilde{b}_K(f)( x_1 \ldots x_{p+1}) =f(x_1\ldots x_{p}) x_{p+1} -(-1)^m x_1 f(x_2 \ldots x_{p+1}),
\end{equation}
\begin{equation} \label{tildekocupA}
(f\tilde{\underset{K}{\smile}} g) (x_1 \ldots x_{p+q}) = (-1)^{mn} f(x_1 \ldots x_p) g(x_{p+1} \ldots  x_{p+q}).
\end{equation}
Let us also define the corresponding cup bracket by
$$[f,g]_{\tilde{\underset{K}{\smile}}}=f\tilde{\underset{K}{\smile}} g -(-1)^{mn} g\tilde{\underset{K}{\smile}} f.$$
\Blm
The product $\tilde{\underset{K}{\smile}}$ is associative and the following formula holds
$$\tilde{b}_K(f)=-[e_A,f]_{\tilde{\underset{K}{\smile}}}$$
for any Koszul cochain $f$ with coefficients in $A$.
\Elm

The proof is immediate. Associativity implies that $[-,-]_{\tilde{\underset{K}{\smile}}}$ is a graded biderivation for the product $\tilde{\underset{K}{\smile}}$. Consequently, one has $\tilde{b}_K(\tilde{b}_K(f))=0$ and
$$\tilde{b}_K(f\tilde{\underset{K}{\smile}} g)=\tilde{b}_K(f) \tilde{\underset{K}{\smile}} g+ (-1)^m f\tilde{\underset{K}{\smile}} \tilde{b}_K(g).$$
Therefore, $(Hom(W_{\bullet},A),\tilde{\underset{K}{\smile}}, \tilde{b}_K )$ is a dga w.r.t. the coefficient weight. The following convention 
is essential for stating the Koszul duality in the next theorem.
\\

\noindent
\emph{Convention}: $(Hom(W_{\bullet},A),\tilde{\underset{K}{\smile}})$ is considered as $\mathbb{N}\times \mathbb{N}$-graded by the \emph{inverse} biweight $(m,p)$.
\\

The homology of the complex $(Hom(W_{\bullet},A), \tilde{b}_K )$ is denoted by $\tilde{HK}^{\bullet}(A)$, it is a unital associative algebra,
 $\mathbb{N}\times \mathbb{N}$-graded by the inverse biweight $(m,p)$. The homogeneous component of biweight $(m,p)$ is denoted by 
 $\tilde{HK}^p(A)_m$. Note that $HK^{\bullet}(A)$ and $\tilde{HK}^{\bullet}(A)$ are different in general. For example, $HK^0(A)=Z(A)$, 
 while $\tilde{HK}^0(A)=\tilde{Z}(A)$ is the graded center of $A$, considering $A$ graded by the weight.
\Bte \label{cohkduality}
Let $V$ be a finite dimensional $k$-vector space, $A=T(V)/(R)$ a quadratic algebra and $A^!=T(V^{\ast})/(R^{\perp})$ the Koszul dual of $A$. 
There is an isomorphism of $\mathbb{N}\times \mathbb{N}$-graded unital associative algebras
\begin{equation} \label{isocohkduality}
(HK^{\bullet}(A), \underset{K}{\smile}) \cong (\tilde{HK}^{\bullet}(A^!), \tilde{\underset{K}{\smile}}).
\end{equation}
In particular, for any $p\geq 0$ and $m\geq 0$, there is a $k$-linear isomorphism
\begin{equation} \label{linearisocohkduality}
HK^p(A)_m \cong \tilde{HK}^m(A^!)_p.
\end{equation}
\Ete
\Bdm
Let us first explain the strategy: it suffices to exhibit a morphism of  $\mathbb{N}\times \mathbb{N}$-graded unital associative algebras
\begin{equation} \label{cochainisocohkduality}
\varphi_A : (Hom(W_{\bullet},A), \underset{K}{\smile}) \rightarrow (Hom(W^!_{\bullet},A^!), \tilde{\underset{K}{\smile}}),
\end{equation}
which is a morphism of complexes w.r.t. $b_K$ and $ \tilde{b}_K$, such that $\varphi_{A^!}\circ \varphi_A=\mathrm{id}$ and 
$\varphi_A\circ \varphi_{A^!}=\mathrm{id}$ -- using the natural isomorphisms $W^{!!}_{\bullet} \cong W_{\bullet}$ and $A^{!!} \cong A$. 
In fact, the isomorphism (\ref{isocohkduality}) will be then given by
$$H(\varphi_A) : (HK^{\bullet}(A), \underset{K}{\smile}) \rightarrow (\tilde{HK}^{\bullet}(A^!), \tilde{\underset{K}{\smile}}).$$

We begin by the definition of $\varphi_A$. Using (\ref{W!p}) and the natural isomorphism  $V^{\ast \otimes p}\cong (V^{\otimes p})^{\ast}$, 
the space $W^!_p$ 
is identified to the orthogonal space of $\sum_{i+2+j=p} V^{\otimes i}\otimes R \otimes V^{\otimes j}$ in $(V^{\otimes p})^{\ast}$. 
The following lemma is standard.
\Blm
For any subspace $F$ of a finite dimensional vector space $E$, denote by $F^{\perp}$ the subspace of $E^{\ast}$ whose elements are the linear 
forms vanishing on $F$. The canonical map $(E/F)^{\ast} \rightarrow E^{\ast}$, transpose of $can: E\rightarrow E/F$, defines an isomorphism
 $(E/F)^{\ast} \cong F^{\perp}$, and the canonical map $E^{\ast} \rightarrow F^{\ast}$, transpose of $can: F\rightarrow E$, defines an 
 isomorphism $E^{\ast}/F^{\perp}  \rightarrow F^{\ast}$.
\Elm
Applying the lemma, we define the $k$-linear isomorphism $\psi_p: W^!_p \rightarrow A_p^{\ast},$
where $A_p^{\ast}$ denotes the dual vector space of
$$A_p=V^{\otimes p} / \sum_{i+2+j=p} V^{\otimes i}\otimes R \otimes V^{\otimes j}.$$
The transpose $\psi_p^{\ast}: A_p \rightarrow  W^{!\ast}_p$ is an isomorphism. Replacing $A$ by $A^!$ and using that $W^{!!}_p \cong W_p$, 
the map $\psi_p^{! \ast}: A^!_p \rightarrow  W^{\ast}_p$ is an isomorphism as well. According to the lemma, $\psi_p^{! \ast}$ is induced by the map sending 
any linear form on $V^{\otimes p}$ to its restriction to $W_p$.
\Bdf \label{varphi}
For any $p\geq 0$, $m\geq 0$ and for any Koszul cochain $f:W_p\rightarrow A_m$, we define the Koszul cochain $\varphi_A(f): W^!_m \rightarrow A^!_p$ 
by the commutative diagram
\begin{eqnarray} \label{diagramvarphi}
W^!_m \ \ & \stackrel{\varphi_A(f)}{\longrightarrow} & A^!_p \nonumber  \\
\downarrow \psi_m &  & \  \downarrow \psi_p^{! \ast} \\
A^{\ast}_m \ \ &  \stackrel{f^{\ast}}{\longrightarrow} & W^{\ast}_p. \nonumber \end{eqnarray}
\Edf

The so-defined $k$-linear map $\varphi_A$ is homogeneous for the biweight $(p,m)$ of $Hom(W_{\bullet},A)$ and the biweight $(m,p)$ 
of $Hom(W^!_{\bullet},A^!)$. Diagram (\ref{diagramvarphi}) applied to $A^!$ and to $\varphi_A(f)$ provides the commutative diagram
\begin{eqnarray} \label{diagramvarphi2}
 W_p\ \  & \stackrel{\varphi_{A^!}(\varphi_A(f))}{\longrightarrow} & A_m \nonumber  \\
\downarrow \psi^!_p &  &  \ \downarrow \psi_m^{\ast} \\
A^{! \ast}_p \ \  &  \stackrel{\varphi_A(f)^{\ast}}{\longrightarrow} & W^{!\ast}_m. \nonumber \end{eqnarray}
Comparing this diagram to the transpose of Diagram (\ref{diagramvarphi}), we obtain $\varphi_{A^!}\circ \varphi_A(f)=f$. 
The proof of $\varphi_A\circ \varphi_{A^!}(h)=h$ for any $h:W^!_m \rightarrow A^!_p$ is similar. So
$$\varphi_A : Hom(W_{\bullet},A) \rightarrow Hom(W^!_{\bullet},A^!)$$
is a $k$-linear isomorphism whose inverse isomorphism is $\varphi_{A^!}$. We continue the proof of Theorem \ref{cohkduality} by the following.
\Bca
The map $\varphi_A $ is an algebra morphism.
\Eca
\Bdm
Let $f:W_p \rightarrow A_m$ and $g:W_q \rightarrow A_n$. For the proof, it is necessary to introduce the cup product \emph{without sign} $\bar{\underset{K}{\smile}}$ defined on $Hom(W_{\bullet},A)$ by
$$(f\bar{\underset{K}{\smile}} g) (x_1 \ldots x_{p+q}) = f(x_1 \ldots x_p) g(x_{p+1} \ldots  x_{p+q}).$$
Conformally to the \emph{ungraded situation} stated in the introduction of this section, the tensor products of linear maps are 
understood \emph{without sign} in the sequel. In particular, the following diagram, whose transpose is used below, commutes.
\begin{eqnarray*}  W_p\otimes W_q  &  \stackrel{f\otimes g}{\longrightarrow} & A_m \otimes A_n \\
 \uparrow \mathrm{can} &  & \ \ \  \downarrow \mu \\
 W_{p+q}\ \  & \stackrel{f\bar{\underset{K}{\smile}} g}{\longrightarrow} & A_{m+n}. \end{eqnarray*}
Tensoring Diagram (\ref{diagramvarphi}) by its analogue for $g$, we write down the commutative diagram
\begin{eqnarray} \label{tensorvarphi}
W^!_m \otimes W^!_n \ \ \ & \stackrel{\varphi_A(f)\otimes \varphi_A(g)}{\longrightarrow} & A^!_p \otimes A^!_q\nonumber  \\
\downarrow \psi_m \otimes \psi_n &  & \ \ \ \downarrow \psi_p^{! \ast} \otimes \psi_q^{! \ast} \\
A^{\ast}_m \otimes A^{\ast}_n \ \ \ &  \stackrel{f^{\ast}\otimes g^{\ast}}{\longrightarrow} & W^{\ast}_p \otimes W^{\ast}_q. \nonumber \end{eqnarray}
Combining this diagram with the following four commutative diagrams
$$\begin{array}{ccc}  W^!_{m+n} \ \ \ & \stackrel{\mathrm{can}}{\longrightarrow} & W^!_m\otimes W^!_n  \\
 \downarrow \psi_{m+n} &  &  \ \ \downarrow \psi_m \otimes \psi_n \\
 A^{\ast}_{m+n}\ \  &  \stackrel{\mu^{\ast}}{\longrightarrow} & A^{\ast}_m \otimes A^{\ast}_n \end{array} \ \ \ \mathrm{and}  \ \ \ \begin{array}{ccc}  A^!_p  \otimes A^!_q \ \ \ \  & \stackrel{\mu^!}{\longrightarrow} & A^!_{p+q}  \\
\downarrow \psi_p^{! \ast} \otimes \psi_q^{! \ast} &  &  \downarrow  \psi_{p+q}^{! \ast}\\
W^{\ast}_p  \otimes W^{\ast}_q \ \ \ \   &  \stackrel{\mathrm{can}}{\longrightarrow} & W^{\ast}_{p+q}
\end{array}
$$
$$\begin{array}{ccc}
W^!_{m+n} & \stackrel{\varphi_A(f) \bar{\underset{K}{\smile}} \varphi_A(g)}{\longrightarrow} & A^!_{p+q}  \\
\downarrow \mathrm{can} &  &  \uparrow \mu^! \\
W^!_m \otimes W^!_n   &  \stackrel{\varphi_A(f) \otimes \varphi_A(g)}{\longrightarrow} & A^!_p \otimes A^!_q \end{array} \ \ \ \mathrm{and}  \ \ \
\begin{array}{ccc} A^{\ast}_m \otimes A^{\ast}_n \ \ \ &  \stackrel{f^{\ast}\otimes g^{\ast}}{\longrightarrow} & W^{\ast}_p \otimes W^{\ast}_q  \\
\uparrow \mu^{\ast} &  &  \downarrow \mathrm{can} \\
A^{\ast}_{m+n}   &  \stackrel{(f\bar{\underset{K}{\smile}} g)^{\ast}}{\longrightarrow} &  W^{\ast}_{p+q}
\end{array}
$$
we obtain the commutativity of
\begin{equation}
\begin{array}{ccc}
W^!_{m+n} & \stackrel{\varphi_A(f) \bar{\underset{K}{\smile}} \varphi_A(g)}{\longrightarrow} & A^!_{p+q}  \\
\downarrow \psi_{m+n} &  &  \downarrow \psi_{p+q}^{! \ast} \\
A^{\ast}_{m+n}   &  \stackrel{(f\bar{\underset{K}{\smile}} g)^{\ast}}{\longrightarrow} &  W^{\ast}_{p+q}.
\end{array}
\end{equation}
Finally, it is sufficient to compare this diagram to Diagram (\ref{diagramvarphi}) applied to $f\bar{\underset{K}{\smile}} g$ instead of $f$, 
for showing that $\varphi_A(f\bar{\underset{K}{\smile}} g)=\varphi_A(f) \bar{\underset{K}{\smile}} \varphi_A(g)$. Multiplying the latter equality 
by $(-1)^{pq}$, we conclude that $\varphi_A(f \underset{K}{\smile} g)=\varphi_A(f) \tilde{\underset{K}{\smile}} \varphi_A(g).$ 
\Edm

Consequently, one has $\varphi_A([f,g]_{\underset{K}{\smile}})=[\varphi_A(f), \varphi_A(g)]_{\tilde{\underset{K}{\smile}}}$.
In particular, $\varphi_A([e_A,f]_{\underset{K}{\smile}})=[e_{A^!}, \varphi_A(f)]_{\tilde{\underset{K}{\smile}}}$, and 
therefore $\varphi_A(b_K(f))=\tilde{b}_K(\varphi_A(f))$ by using the fundamental formulas. Theorem \ref{cohkduality} is thus proved.
\Edm
\\

We illustrate Theorem \ref{cohkduality} by the example $A=k[x]$, that is $V=k.x$ and $R=0$. 
The Koszul dual of $A$ is $A^!=k\oplus k.x^{\ast}$ with $x^{\ast 2}=0$. It is straightforward to verify the following isomorphisms for any $m\geq 0$
$$HK^0(A)_m \cong k.(1\mapsto x^m) \cong k.(x^{\ast m} \mapsto 1) \cong \tilde{HK}^m(A^!)_0$$
$$HK^1(A)_m \cong k.(x\mapsto x^m) \cong k.(x^{\ast m} \mapsto x^{\ast}) \cong \tilde{HK}^m(A^!)_1$$
$$HK^p(A)_m \cong 0 \cong \tilde{HK}^m(A^!)_p \ \mathrm{for \ any} \ p\geq 2,$$
and it is also direct to check that the products work well. Remark that $HK^0(A)_m$ is not isomorphic to $\tilde{HK}^0(A^!)_m$ for any $m\geq 2$, 
so the exchange $p \leftrightarrow m$ is essential in Theorem \ref{cohkduality}. Passing to the modified version $\tilde{HK}^m(A^!)_p$ is also essential, since $HK^m(A^!)_0$ is not isomorphic to $HK^0(A)_m$ when $m$ is odd. Moreover, it is clear that $HK^0(A^!) \ncong HK^0(A)$.

\subsection{Koszul duality in higher cohomology}

As in Subsection \ref{subsechkoco}, we define the tilde version of the Koszul higher cohomology. Clearly,
 $e_A\tilde{\underset{K}{\smile}}e_A=0$, so that $e_A\tilde{\underset{K}{\smile}} -$ is a cochain differential on $Hom(W_{\bullet},A)$. 
 Next, $\overline{e}_A\tilde{\underset{K}{\smile}} -$ is a cochain differential on $\tilde{HK}^{\bullet}(A)$ denoted 
 by $\tilde{\partial}_{\smile}$. The homology of $\tilde{HK}^{\bullet}(A)$ endowed with $\tilde{\partial}_{\smile}$ is denoted 
 by $\tilde{HK}_{hi}^{\bullet}(A)$. The associative algebra $(\tilde{HK}_{hi}^{\bullet}(A), \tilde{\underset{K}{\smile}})$ is $\mathbb{N}\times \mathbb{N}$-graded by the inverse biweight. Since 
 $$H(\varphi_A)(\overline{e}_A \underset{K}{\smile} \alpha)=\overline{e}_{A^!} \tilde{\underset{K}{\smile}} H(\varphi_A)(\alpha),$$
for any $\alpha$ in $HK^{\bullet}(A)$, Theorem \ref{cohkduality} implies that the isomorphism 
$H(\varphi_A) : HK^{\bullet}(A) \rightarrow \tilde{HK}^{\bullet}(A^!)$ is also an isomorphism of complexes w.r.t. the differentials 
$\partial_{\smile}$ and $\tilde{\partial}_{\smile}$. We have thus proved the following higher Koszul duality theorem.
\Bte \label{hcohkduality}
Let $V$ be a finite dimensional $k$-vector space and $A=T(V)/(R)$ a quadratic algebra. Let $A^!=T(V^{\ast})/(R^{\perp})$ be the Koszul dual of $A$.
 There is an isomorphism of $\mathbb{N}\times \mathbb{N}$-graded associative algebras
\begin{equation} \label{isohcohkduality}
(HK_{hi}^{\bullet}(A), \underset{K}{\smile}) \cong (\tilde{HK}_{hi}^{\bullet}(A^!), \tilde{\underset{K}{\smile}}).
\end{equation}
In particular, for any $p\geq 0$ and $m\geq 0$, there is a $k$-linear isomorphism
\begin{equation} \label{linearisocohkduality}
HK_{hi}^p(A)_m \cong \tilde{HK}_{hi}^m(A^!)_p.
\end{equation}
\Ete

\subsection{Koszul duality in homology}

We proceed as we have done for cohomology in Subsection \ref{subseckdco}. We define a modified version of 
Koszul homology by exchanging homological and 
coefficient weights. Precisely, for $f:W_p \rightarrow A_m$ and 
$z= a\otimes x_1 \ldots x_{q}$ in $ A_n\otimes W_q$, we define $\tilde{b}_K(z)$, $f \tilde{\underset{K}{\frown}}z$ and $z\tilde{\underset{K}{\frown}}f$ 
by
\begin{equation} \label{tildedefb}
\tilde{b}_K(z) = ax_1\otimes x_2 \ldots x_q +(-1)^n x_q a \otimes x_1 \ldots x_{q-1},
\end{equation}
\begin{equation} \label{tildekolcapA}
  f \tilde{\underset{K}{\frown}} z = (-1)^{(n-m)m} f(x_{q-p+1} \ldots x_q)a \otimes x_1 \ldots  x_{q-p},
\end{equation}
\begin{equation} \label{tildekorcapA}
  z \tilde{\underset{K}{\frown}} f = (-1)^{mn} af(x_1 \ldots x_p) \otimes x_{p+1} \ldots  x_q.
\end{equation}
The corresponding cap bracket is
$$[f,z]_{\tilde{\underset{K}{\frown}}}=f\tilde{\underset{K}{\frown}} z -(-1)^{mn} z\tilde{\underset{K}{\frown}} f.$$
It is just routine to verify the following associativity relations:
\begin{eqnarray*}
  f\tilde{\underset{K}{\frown}}\, (g\tilde{\underset{K}{\frown}} z) & = & (f\tilde{\underset{K}{\smile}}g)\, \tilde{\underset{K}{\frown}} z,\\
  (z\tilde{\underset{K}{\frown}} g)\, \tilde{\underset{K}{\frown}} f & = & z \tilde{\underset{K}{\frown}}\, (g\tilde{\underset{K}{\smile}} f),\\
  f\tilde{\underset{K}{\frown}}\, (z\tilde{\underset{K}{\frown}} g) & = & (\tilde{f\underset{K}{\frown}}z)\, \tilde{\underset{K}{\frown}} g,
\end{eqnarray*}
and the fundamental formula $$\tilde{b}_K(z)=-[e_A,z]_{\tilde{\underset{K}{\frown}}}.$$

The associativity relations imply that $[-,-]_{\tilde{\underset{K}{\frown}}}$ is a graded biderivation for the product 
$\tilde{\underset{K}{\smile}}$ in the first argument and the actions $\tilde{\underset{K}{\frown}}$ in the second argument. 
From that, it is straightforward to deduce $\tilde{b}_K(\tilde{b}_K(z))=0$ and
$$\tilde{b}_K(f\tilde{\underset{K}{\frown}} z)=\tilde{b}_K(f)\, \tilde{\underset{K}{\frown}} z + (-1)^m f\tilde{\underset{K}{\frown}}\, \tilde{b}_K(z),$$
$$\tilde{b}_K(z\tilde{\underset{K}{\frown}} f)=\tilde{b}_K(z)\, \tilde{\underset{K}{\frown}} f + (-1)^n z\tilde{\underset{K}{\frown}}\, \tilde{b}_K(f).$$
Therefore, $(A \otimes W_{\bullet},\tilde{\underset{K}{\frown}}, \tilde{b}_K )$ is a differential graded bimodule w.r.t. the coefficient weight 
over the dga $(Hom(W_{\bullet},A), \tilde{\underset{K}{\smile}}, \tilde{b}_K )$.

The homology of the complex $(A \otimes W_{\bullet}, \tilde{b}_K )$ is denoted by $\tilde{HK}_{\bullet}(A)$. It is a
 $\tilde{HK}^{\bullet}(A)$-bimodule, $\mathbb{N}\times \mathbb{N}$-graded by the \emph{inverse} biweight. The homogeneous component of biweight 
 $(n,q)$ is denoted by $\tilde{HK}_q(A)_n$. Note that $HK_0(A)_0 \cong k$, while $\tilde{HK}_0(A)_0 \cong 0$ if char$(k)\neq 2$.

In order to state the Koszul duality in homology, we need to slightly generalize the formalism described up to now in this section, by replacing 
the graded space of coefficients, namely $A$, by an arbitrary $\mathbb{Z}$-graded $A$-bimodule $M$,  whose degree is still called the \emph{weight}.
The formalism described up to now for $M=A$ extends immediately to such a graded $M$ by using \emph{the same} $b_K$, $\underset{K}{\smile}$,
 $\underset{K}{\frown}$, $\tilde{b}_K$, $\tilde{\underset{K}{\smile}}$, $\tilde{\underset{K}{\frown}}$. We obtain the following general formalism.
\begin{enumerate}
\item $Hom(W_{\bullet},M)$ is a $(Hom(W_{\bullet},A),\underset{K}{\smile})$-bimodule for $\underset{K}{\smile}$, $\mathbb{N}\times \mathbb{Z}$-graded 
by the biweight, and $HK^{\bullet}(A,M)$ is a $\mathbb{N}\times \mathbb{Z}$-graded $(HK^{\bullet}(A),\underset{K}{\smile})$-bimodule.
\item $Hom(W_{\bullet},M)$ is a $(Hom(W_{\bullet},A),\tilde{\underset{K}{\smile}})$-bimodule for $\tilde{\underset{K}{\smile}}$, 
$\mathbb{Z}\times \mathbb{N}$-graded by the inverse biweight, and $\tilde{HK}^{\bullet}(A,M)$ is a $\mathbb{Z}\times \mathbb{N}$-graded $(\tilde{HK}^{\bullet}(A),\tilde{\underset{K}{\smile}})$-bimodule.
\item $M\otimes W_{\bullet}$ is a $(Hom(W_{\bullet},A),\underset{K}{\smile})$-bimodule for $\underset{K}{\frown}$, 
$\mathbb{N}\times \mathbb{Z}$-graded by the biweight, and $HK_{\bullet}(A,M)$ is a $\mathbb{N}\times \mathbb{Z}$-graded
 $(HK^{\bullet}(A),\underset{K}{\smile})$-bimodule.
\item $M\otimes W_{\bullet}$ is a $(Hom(W_{\bullet},A),\tilde{\underset{K}{\smile}})$-bimodule for $\tilde{\underset{K}{\frown}}$, 
$\mathbb{Z}\times \mathbb{N}$-graded by the inverse biweight, and $\tilde{HK}_{\bullet}(A,M)$ is a $\mathbb{Z}\times \mathbb{N}$-graded
 $(\tilde{HK}^{\bullet}(A),\tilde{\underset{K}{\smile}})$-bimodule.
\end{enumerate}

Apart from the case $M=A$, we will need to consider the graded dual $M=A^{\ast}= \bigoplus_{m\geq 0} A^{\ast}_m$. It would be more natural to 
grade $A^{\ast}$ by the weight $-m$, but in order to avoid notational complications, we prefer to use the nonnegative weight $m$. So all the 
biweights used below will belong to $\mathbb{N}\times \mathbb{N}$. We recall the actions of the graded $A$-bimodule $A^{\ast}$. 
For any $u$ in $A^{\ast}_m$ and $a$ in $A_n$, they are defined by $a.u$ and $u.a$ in $A^{\ast}_{m-n}$, where  
\begin{equation} \label{leftactionAdual}
(a.u)(a') = (-1)^n u(a'a),
\end{equation}
\begin{equation} \label{rightactionAdual}
(u.a)(a') = u(aa'),
\end{equation}
for any $a'$ in $A_{m-n}$. We are now 
ready to state the following Koszul duality theorem in homology, completing Theorem \ref{cohkduality}.
\Bte \label{hkduality}
Let $V$ be a finite dimensional $k$-vector space and $A=T(V)/(R)$ a quadratic algebra. Let $A^!=T(V^{\ast})/(R^{\perp})$ be the Koszul dual of $A$. 
There is an isomorphism \begin{equation} \label{isohkduality}
HK_{\bullet}(A) \cong \tilde{HK}^{\bullet}(A^!, A^{!\ast}), \end{equation}
from the $(HK^{\bullet}(A),\underset{K}{\smile})$-bimodule $HK_{\bullet}(A)$ with actions $\underset{K}{\frown}$, 
$\mathbb{N}\times \mathbb{N}$-graded by the biweight, to the $(\tilde{HK}^{\bullet}(A^!),\tilde{\underset{K}{\smile}})$-bimodule
 $\tilde{HK}^{\bullet}(A^!, A^{!\ast})$ with actions $\tilde{\underset{K}{\smile}}$, $\mathbb{N}\times \mathbb{N}$-graded by the inverse biweight.  
 In particular, for any $p\geq 0$ and $m\geq 0$, there is a $k$-linear isomorphism
\begin{equation} \label{linearisohkduality}
HK_p(A)_m \cong \tilde{HK}^m(A^!, A^{!\ast})_p.
\end{equation}
\Ete
\Bdm
It is sufficient to exhibit an isomorphism \begin{equation} \label{chainisohkduality}
\theta_A : A \otimes W_{\bullet} \rightarrow Hom(W^!_{\bullet},A^{!\ast}),
\end{equation}
from the $(Hom(W_{\bullet},A),\underset{K}{\smile})$-bimodule $A \otimes W_{\bullet}$ with actions $\underset{K}{\frown}$ to 
the $(Hom(W^!_{\bullet},A^!),\underset{K}{\smile})$-bimodule $Hom(W^!_{\bullet},A^{!\ast})$ with actions $\tilde{\underset{K}{\smile}}$, 
such that $\theta_A$ is homogeneous for the biweights as in the statement and $\theta_A$ is a morphism of complexes w.r.t. $b_K$ and 
$ \tilde{b}_K$. After doing so, the isomorphism (\ref{isohkduality}) will be given by
$$H(\theta_A) : HK_{\bullet}(A) \cong \tilde{HK}^{\bullet}(A^!, A^{!\ast}).$$ 
For defining the linear map $\theta_A : A_m \otimes W_p \rightarrow Hom(W^!_m,A^{!\ast}_p)$, we use the linear isomorphisms  
$\psi_m^{\ast}: A_m \rightarrow  W^{!\ast}_m$ and  $\psi_p^!: W_p \rightarrow  A^{!\ast}_p$ defined in the proof of Theorem \ref{cohkduality}. 
For any $z= a\otimes x_1 \ldots x_p$ in $ A_m\otimes W_p$, set
\begin{equation} \label{explicitchainisohkduality}
\theta_A(z)(w)= \psi_m^{\ast}(a)(w)\, \psi_p^!(x_1 \ldots x_p),
\end{equation}
for any $w$ in $W^!_m$. The so-defined linear map $\theta_A$ is homogeneous for the biweight of $A \otimes W_{\bullet}$ and the inverse biweight 
of $Hom(W^!_{\bullet},A^{!\ast})$.

Defining
\begin{equation} \label{inversechainisohkduality}
\theta'_A :  Hom(W^!_m,A^{!\ast}_p) \rightarrow A_m \otimes W_p,
\end{equation}
by $\theta'_A(f)= \sum_{i \in I} e_i \otimes (\psi_p^{!-1}\circ f \circ \psi_m^{-1}(e_i^{\ast}))$ for any linear 
$f:W^!_m \rightarrow A^{!\ast}_p$, where $(e_i)_{i \in I}$ is a basis of the space $A_m$ and $(e_i^{\ast})_{i \in I}$ is its dual basis, it is easy to verify that $\theta_A$ is an 
isomorphism whose inverse is $\theta'_A$.
\Bca
Using $\varphi_A$, consider the $Hom(W^!_{\bullet},A^!)$-bimodule $Hom(W^!_{\bullet},A^{!\ast})$ as a $Hom(W_{\bullet},A)$-bimodule. 
The map $\theta_A: A \otimes W_{\bullet} \rightarrow Hom(W^!_{\bullet},A^{!\ast})$ is a morphism of $Hom(W_{\bullet},A)$-bimodules.
\Eca
\Bdm
This amounts to prove that
\begin{equation} \label{thetalaction}
\theta_A(f\underset{K}{\frown}z) = \varphi_A(f) \tilde{\underset{K}{\smile}} \theta_A(z),
\end{equation}
\begin{equation} \label{thetaraction}
\theta_A(z\underset{K}{\frown}f) = \theta_A(z) \tilde{\underset{K}{\smile}} \varphi_A(f),
\end{equation}
for any $z= a\otimes x_1 \ldots x_p$ in $ A_m\otimes W_p$ and $f:W_q \rightarrow A_n$, with $p\geq q$.

Analogously to $\bar{\underset{K}{\smile}}$, define the cap products without sign $\bar{\underset{K}{\frown}}$. Firstly we prove 
\begin{equation} \label{thetalbaraction}
\theta_A(f \bar{\underset{K}{\frown}} z) = \varphi_A(f) \bar{\underset{K}{\smile}} \theta_A(z),
\end{equation}
leaving to the reader the proof of
\begin{equation} \label{thetarbaraction}
\theta_A(z \bar{\underset{K}{\frown}} f) = \theta_A(z) \bar{\underset{K}{\smile}} \varphi_A(f).
\end{equation}
For any $w=y_1 \ldots y_{m+n} \in W^!_{m+n}$, we deduce from equality (\ref{explicitchainisohkduality}) that
$$\theta_A(f \bar{\underset{K}{\frown}} z)(w)= \psi_{m+n}^{\ast}(f(x_{p-q+1} \ldots x_p)a)(w)\, \psi_{p-q}^!(x_1 \ldots x_{p-q}).$$
Write $w=w_1w_2$ where $w_1=y_1 \ldots y_n \in W^!_n$ and $w_2=y_{n+1} \ldots y_{m+n} \in W^!_m$, so that
$$\theta_A(z)(w_2)= \psi_m^{\ast}(a)(w_2)\, \psi_p^!(x_1 \ldots x_p).$$
Denoting by $\bar{.}$ the left action of an element of $A^!_q$ on an element $A^{!\ast}_p$ giving an element of $A^{!\ast}_{p-q}$ as 
in (\ref{leftactionAdual}) \emph{but without sign}, we have
\begin{eqnarray*}
(\varphi_A(f) \bar{\underset{K}{\smile}}\theta_A(z)) (w) & = & \varphi_A(f)(w_1)\, \bar{.}\, \theta_A(z)(w_2) \\
                             & = & \psi_q^{!\ast -1}\circ f^{\ast} \circ \psi_n(w_1)\, \bar{.}\, (\psi_m^{\ast}(a)(w_2)\, \psi_p^!(x_1 \ldots x_p))\\
                             & = & \psi_m^{\ast}(a)(w_2)\, (\psi_q^{!\ast -1}\circ f^{\ast} \circ \psi_n(w_1)\, \bar{.}\,  \psi_p^!(x_1 \ldots x_p)).
\end{eqnarray*}

Next, for any $a' \in A^!_{p-q}$, one has
$$(\psi_q^{!\ast -1}\circ f^{\ast} \circ \psi_n(w_1)\,\bar{.}\,\psi_p^!(x_1 \ldots x_p))(a') = \psi_p^!(x_1 \ldots x_p)(a'(\psi_q^{!\ast -1}\circ f^{\ast} \circ \psi_n(w_1))).$$
The right-hand side is equal to
 $\psi_{p-q}^!(x_1 \ldots x_{p-q})(a')\, \psi_q^!(x_{p-q+1} \ldots x_p)(\psi_q^{!\ast -1}\circ f^{\ast} \circ \psi_n(w_1))$, 
 by using the commutative diagram 
 $$\begin{array}{ccc}  W_p \ \ \ & \stackrel{\mathrm{can}}{\longrightarrow} & W_{p-q}\otimes W_q  \\
 \downarrow \psi^!_p &  &  \ \ \downarrow \psi^!_{p-q} \otimes \psi^!_q \\
 A^{! \ast}_p\ \  &  \stackrel{\mu^{! \ast}}{\longrightarrow} & A^{! \ast}_{p-q} \otimes A^{! \ast}_q.
\end{array} $$
Therefore, we obtain
$$(\varphi_A(f) \bar{\underset{K}{\smile}}\theta_A(z)) (w)=\psi_m^{\ast}(a)(w_2)\,\psi_q^!(x_{p-q+1} \ldots x_p)(\psi_q^{!\ast -1}\circ f^{\ast} 
\circ \psi_n(w_1))\,\psi_{p-q}^!(x_1 \ldots x_{p-q}).$$
By duality, $\psi_q^!(x_{p-q+1} \ldots x_p)(\psi_q^{!\ast -1}\circ f^{\ast} \circ \psi_n(w_1))$ is equal to 
$\psi^{\ast}_n(f(x_{p-q+1} \ldots x_p))(w_1)$. Moreover, the commutative diagram 
$$\begin{array}{ccc}  A_n  \otimes A_m \ \ \ \  & \stackrel{\mu}{\longrightarrow} & A_{m+n}  \\
\downarrow \psi_n^{\ast} \otimes \psi_m^{\ast} &  &  \downarrow  \psi_{m+n}^{\ast}\\
W^{! \ast}_n  \otimes W^{! \ast}_m \ \ \ \   &  \stackrel{\mathrm{can}}{\longrightarrow} & W^{! \ast}_{m+n}
\end{array} $$
shows that
$$\psi^{\ast}_n(f(x_{p-q+1} \ldots x_p))(w_1)\,\psi_m^{\ast}(a)(w_2)=\psi^{\ast}_{m+n} (f(x_{p-q+1} \ldots x_p)a)(w_1w_2).$$
Thus equality (\ref{thetalbaraction}) is proved. We draw the following
$$\theta_A(f \underset{K}{\frown} z) =  (-1)^{pq} (-1)^q \varphi_A(f) \bar{\underset{K}{\smile}}\,\theta_A(z).$$
Recall that $\varphi_A(f): W^!_n \rightarrow A^!_q$ and $\theta_A(z):W^!_m \rightarrow A^{! \ast}_p$, so that $(-1)^{pq}$ is equal to the sign 
defining $\tilde{\underset{K}{\smile}}$ from $\bar{\underset{K}{\smile}}$, without forgetting the sign $(-1)^q$ defining the left action of $A^!_q$ 
on $A^{! \ast}_p$ as in (\ref{leftactionAdual}). 
Therefore $\theta_A(f \underset{K}{\frown} z) =  \varphi_A(f) \tilde{\underset{K}{\smile}} \theta_A(z)$.

Similarly, 
$\theta_A(z \underset{K}{\frown} f)  =  (-1)^{pq}\, \theta_A(z) \bar{\underset{K}{\smile}}\, \varphi_A(f) =  \theta_A(z) \tilde{\underset{K}{\smile}}\, \varphi_A(f)$.
\Edm

\medskip

Consequently, one gets $\theta_A([f,z]_{\underset{K}{\frown}})=[\varphi_A(f), \theta_A(z)]_{\tilde{\underset{K}{\smile}}}$, and $\theta_A(b_K(z))=\tilde{b}_K(\theta_A(z))$ by using the fundamental formulas. Theorem \ref{hkduality} is proved.
\Edm

\Brm 
\Erm
Denote by $\mathcal{C}$ the Manin category of quadratic $k$-algebras over finite dimensional vector spaces, and by $\mathcal{E}$ the category 
of the $\mathbb{N}\times \mathbb{N}$-graded $k$-vector spaces whose components are finite dimensional. We know that $A\mapsto HK_{\bullet}(A)$ is a covariant functor $F$ from $\mathcal{C}$ to $\mathcal{E}$. Moreover, $A\mapsto HK^{\bullet}(A, A^{\ast})$
is a contravariant functor $G$ from $\mathcal{C}$ to $\mathcal{E}$ where $A^{\ast}$ is the graded dual, hence the same holds for 
$\tilde{G}: A\mapsto \tilde{HK}^{\bullet}(A, A^{\ast})$. The proof of Theorem \ref{hkduality} shows that the duality functor $D: A \mapsto A^!$ 
defines a \emph{natural isomorphism} $\theta$ from $F$ to $\tilde{G} \circ D$.

\subsection{Koszul duality in higher homology}

Generalizing the modified version of higher Koszul cohomology to any $\mathbb{Z}$-graded bimodule $M$, we obtain the following higher Koszul duality theorem in homology, completing Theorem \ref{hcohkduality}.
\Bte \label{hhkduality}
Let $V$ be a finite dimensional $k$-vector space and $A=T(V)/(R)$ a quadratic algebra. Let $A^!=T(V^{\ast})/(R^{\perp})$ be the Koszul dual of $A$.
 There is an isomorphism of $\mathbb{N}\times \mathbb{N}$-graded $HK_{hi}^{\bullet}(A)$-bimodules \begin{equation} \label{isohhkduality}
 HK^{hi}_{\bullet}(A) \cong \tilde{HK}_{hi}^{\bullet}(A^!, A^{!\ast}).
\end{equation}
In particular, for any $p\geq 0$ and $m\geq 0$, there is a $k$-linear isomorphism
\begin{equation} \label{linearisocohkduality}
HK^{hi}_p(A)_m \cong \tilde{HK}_{hi}^m(A^!, A^{!\ast})_p.
\end{equation}
\Ete

\setcounter{equation}{0}

\section{A non-Koszul example} \label{nkexample}

\subsection{Koszul algebras with two generators}

The Koszul algebras with two generators were explicitly determined by the first author in~\cite{rb:wcqa}. The result is recalled below without 
proof. The paper~\cite{rb:wcqa} was devoted to study changes of generators in quadratic algebras and their consequences on confluence. The result 
was obtained by using Priddy's theorem, which asserts that any weakly confluent quadratic algebra is Koszul, and some lattice techniques for the converse ``Koszulity implies strong confluence'' in case of two generators and two relations.

Assume that $V=k.x\oplus k.y$, $R$ is a subspace of $V\otimes V$ and $A=T(V)/(R)$. If $R=0$ or $R=V\otimes V$, then $A$ is Koszul. 
If $\dim(R)=1$, then $A$ is Koszul according to Gerasimov's theorem~\cite{gera:distrib, rb:gera}. If $\dim(R)=3$, $A$ is Koszul since $\dim(R^{\perp})=1$ and $A^!$ is Koszul. For two relations, the 
Koszul algebras are given by the following proposition.
\Bpo
Under the previous assumptions and identifying $A$ to its quadratic relations, the Koszul algebras with two generators and two relations are 
the following.
\begin{equation}
\left\{ \begin{array}{lll}
xy & = & 0\\ x^2 & = & 0
\end{array}
\right.
\ \mbox{and}\ \left\{ \begin{array}{lll}
yx & = & \alpha xy\\ x^2 & = & 0
\end{array}
\right.
\mbox{are\ Koszul}.
\end{equation}
\begin{equation}
\left\{ \begin{array}{lll}
yx & = & \alpha x^2 \\ xy & = & \beta x^2
\end{array}
\right.
\ \mbox{is \ Koszul}\ \Leftrightarrow \ \alpha = \beta.
\end{equation}
\begin{equation} \label{examplesfamily}
\left\{ \begin{array}{lll}
y^2 & = & \alpha xy + \beta yx \\ x^2 & = & 0
\end{array}
\right.
\ \mbox{is \ Koszul}\ \Leftrightarrow \ \alpha = \beta.
\end{equation}
\begin{equation}
\left\{ \begin{array}{lll}
y^2 & = & \alpha x^2 + \beta yx\\ xy & = & \gamma x^2
\end{array}
\right.
\ \mbox{is \ Koszul}\ \Leftrightarrow \ \alpha = 0 \ \mbox{and}\ \beta=\gamma.
\end{equation}
\begin{equation}
\left\{ \begin{array}{lll}
y^2 & = & \alpha x^2 + \beta xy\\ yx & = & \gamma x^2 + \delta xy
\end{array}
\right.
\ \mbox{is \ Koszul}\ \Leftrightarrow \ \left\{ \begin{array}{lll}
\beta (1-\delta) & = & \gamma (1+ \delta) \\ \alpha (1-\delta ^2) & = & -\beta \gamma \delta
\end{array}
\right.
\end{equation}
\Epo

Throughout the remainder of this section, $A$ denotes the non-Koszul quadratic algebra
$$A=k\langle x,y \rangle /\langle x^2, y^2-xy\rangle .$$
It is immediate that the cubic relations $y^3=xy^2=yxy=0$ and $y^2x=xyx$ hold in $A$. Moreover, $A_3$ is $1$-dimensional generated by $xyx$ and $A_m=0$ for any $m\geq 4$. Therefore, $\dim(A)=6$ and $1$, $x$, $y$, $xy$, $yx$, $xyx$ form a linear basis of $A$. This basis will be continually used during the rather long but routine calculations of the various homology and cohomology spaces. We will just state the results, assuming that the characteristic of $k$ is zero. It is easy to show that $W_p=k.x^p$ for any $p\geq 3$.

\subsection{The Koszul homology of $A$}

The complex of Koszul chains of $A$ with coefficients in $A$ is given by
\begin{equation} \label{hkcomplexexample}
\ldots \stackrel{b_K}\longrightarrow A\otimes x^4 \stackrel{b_K}\longrightarrow
A\otimes  x^3 \stackrel{b_K}\longrightarrow A\otimes R \stackrel{b_K}\longrightarrow
A\otimes  V \stackrel{b_K}\longrightarrow A \longrightarrow 0,
\end{equation} where the maps $b_K$ are successively given by
$$b_K(a\otimes x)=ax-xa \ \mbox{and} \ b_K(a'\otimes y)=a'y-ya',$$
$$b_K(a\otimes x^2)= (ax + xa)\otimes x \ \mbox{and}\ b_K(a'\otimes(y^2-xy))=-ya'\otimes x + (a'y+ya'-a'x)\otimes y,$$
$$b_K(a\otimes x^p)=(ax+(-1)^p xa)\otimes x^{p-1},$$
for any $a$, $a'$ in $A$, and $p\geq 3$. 
\Bpo   \label{hkexample}
The Koszul homology of $A$ is given by
\begin{enumerate}
\item $HK_0(A)$ is $4$-dimensional, generated by the classes of $1$, $x$, $y$ and $xy$,

\item $HK_1(A)$ is $3$-dimensional, generated by the classes of $1\otimes x$, $1 \otimes y$ and $y\otimes y$,

\item $HK_2(A)$ is $3$-dimensional, generated by the classes of $x\otimes x^2$, $yx\otimes x^2 + (xy+yx) \otimes (y^2-xy)$ and $xyx\otimes (y^2-xy)$,

\item for any $p\geq 3$ odd (resp. even), $HK_p(A)$ is $1$-dimensional, generated by the class of $1\otimes x^p$ (resp. $x\otimes x^p$).
\end{enumerate}
\Epo
\Bpo  \label{hkhexample}
The higher Koszul homology of $A$ is given by
\begin{enumerate}
\item $HK^{hi}_0(A)\cong k$,
\item $HK^{hi}_1(A)\cong 0$,
\item $HK^{hi}_2(A)$ is $2$-dimensional, generated by the classes of $[yx\otimes x^2 + (xy+yx) \otimes (y^2-xy)]$ and $[xyx\otimes (y^2-xy)]$,
\item $HK^{hi}_p(A)\cong 0$ for any $p\geq 3$.
\end{enumerate}
\Epo

The next proposition shows that $A$ satisfies Conjecture \ref{lconhkoszul}.
\Bpo
The homology of the complex $K_{\ell}(A)$ is given by
\begin{enumerate}
\item $H_0(K_{\ell}(A))\cong k$,
\item $H_1(K_{\ell}(A))\cong 0$,
\item $H_2(K_{\ell}(A))$ is $2$-dimensional, generated by the classes of $yx\otimes (y^2-xy)$ and $xyx\otimes (y^2-xy)$,
\item $H_p(K_{\ell}(A))\cong 0$ for any $p\geq 3$.
\end{enumerate}
\Epo

\subsection{The Koszul cohomology of $A$}

Recall that for any finite dimensional vector space $E$, the linear map $can: A\otimes E^{\ast} \rightarrow Hom(E,A)$ defined by
$can (a \otimes u)(x)=u(x)a$ for any $a$ in $A$, $u$ in $E^{\ast}$ and $x$ in $E$, is an isomorphism. Using this, define the isomorphism of complexes
$$can: A \otimes W_{\bullet}^{\ast} \cong Hom(W_{\bullet},A).$$
The differential of $A \otimes W_{\bullet}^{\ast}$ is obtained by carrying the differential $b_K$ of $Hom(W_{\bullet},A)$, and is still 
denoted by $b_K$.

The dual basis of $V^{\ast}$ corresponding to the basis $(x,y)$ of $V$ is $(x^{\ast},y^{\ast})$. Denote by $x^{\ast 2}$ the restriction to $R$ 
of the linear form $x^{\ast} \otimes x^{\ast}$ on $V\otimes V$, and analogously for $x^{\ast}y^{\ast}$, $y^{\ast}x^{\ast}$ and $y^{\ast 2}$. 
Clearly $x^{\ast 2}$ and $y^{\ast 2}$ form a basis of $R^{\ast}$, and we have the following relations in $R^{\ast}$:
$$x^{\ast}y^{\ast}=-y^{\ast 2}, \ y^{\ast}x^{\ast}=0.$$
For any $p\geq 3$, denote by $x^{\ast p}$ the restriction to $W_p$ of the linear form $x^{\ast\otimes p}$ on $V^{\otimes p}$, so that $W_p^{\ast}$ 
is generated by $x^{\ast p}$. Then it is routine to write down the complex $(A \otimes W_{\bullet}^{\ast}, b_K)$, and to get the following. 
\Bpo \label{cohkexample}
The Koszul cohomology of $A$ is given by
\begin{enumerate}
\item $HK^0(A)$ is $2$-dimensional, generated by $1$ and $xyx$,
\item $HK^1(A)$ is $2$-dimensional, generated by the classes of $x\otimes x^{\ast}+y\otimes y^{\ast} \cong e_A$ and $xy\otimes y^{\ast}$,
\item $HK^2(A)$ is $4$-dimensional, generated by the classes of $1\otimes x^{\ast 2}$, $1\otimes y^{\ast 2}$, $y\otimes y^{\ast 2}$ and 
$xyx \otimes y^{\ast 2}$,
\item for any $p\geq 3$ odd (resp. even), $HK^p(A)$ is $1$-dimensional, generated by the class of $x\otimes x^{\ast p}$ (resp. $1\otimes x^{\ast p}$).
\end{enumerate}\Epo
\Bpo  \label{cohkhexample}
The higher Koszul cohomology of $A$ is given by
\begin{enumerate}
\item $HK_{hi}^0(A)$ is $1$-dimensional, generated by the class of $xyx$,
\item $HK_{hi}^1(A)$ is $1$-dimensional, generated by the class of $[xy\otimes y^{\ast}]$,
\item $HK_{hi}^2(A)$ is $3$-dimensional, generated by the classes of $[1\otimes y^{\ast 2}]$, $[y\otimes y^{\ast 2}]$ and $[xyx\otimes y^{\ast 2}]$,
\item $HK_{hi}^p(A)\cong 0$ for any $p\geq 3$.
\end{enumerate}
\Epo

We do not know whether the following proposition holds or not for any quadratic algebra.
\Bpo
The algebra $(HK^{\bullet}(A), \underset{K}{\smile})$ is graded commutative. The $(HK^{\bullet}(A), \underset{K}{\smile})$-bimodule 
$HK_{\bullet}(A)$ is graded symmetric for the actions $\underset{K}{\frown}$. 
\Epo

We leave the verifications of this proposition to the reader by calculating the cup and cap products of the explicit classes given in Proposition \ref{hkexample} and Proposition \ref{cohkexample}. In higher Koszul cohomology, the products of two biweight homogeneous classes vanish, except
$$[xyx]\underset{K}{\smile}[[1\otimes y^{\ast 2}]]=[[1\otimes y^{\ast 2}]]\underset{K}{\smile}[xyx]=[[xyx\otimes y^{\ast 2}]].$$

Examining the possible biweights, we see also that the higher Koszul cohomology of $A$ acts on the higher Koszul homology of $A$ by zero.

\subsection{The Hochschild (co)homology of $A$}

Apart from standard examples including Koszul algebras, it is difficult to compute the Hochschild (co)homology of an associative algebra given 
by generators and relations. The bar resolution is too large and, if the algebra is graded, a construction of the minimal projective resolution 
is too hard to perform in general. Fortunately, in case of monomial relations, Bardzell's resolution provides a minimal projective resolution 
whose calculation is tractable. The differential and the contracting homotopy of Bardzell's resolution are simultaneously defined in homological 
degree $p$ from ($p-1$)-ambiguities. The ambiguities are monomials simply defined from the well-chosen reduction system $\mathcal{R}$ defining 
the algebra.

The third author and Chouhy have extended Bardzell's resolution to any algebra, not necessarily graded, defined by relations on a finite
 quiver~\cite{cs:reso}. Guiraud, Hoffbeck and Malbos~\cite{ghm:polygraph} have constructed a resolution which may be related to 
 the construction of~\cite{cs:reso}. The first step consists in well-choosing a reduction system $\mathcal{R}$ of the algebra $A$. The resolution $S(A)$ of~\cite{cs:reso} is in some sense 
 a deformation of Bardzell's one. The bimodules of the resolution $S(A)$ are free, and the free bimodule in homological degree $p$ is generated by 
 the ($p-1$)-ambiguities of the associated monomial algebra. The differential and the contracting homotopy are simultaneously defined by induction 
 on $p$. We apply this construction to our favorite non-Koszul algebra $A$, without giving the details.  

The construction of $S(A)$ starts with $x<y$, the corresponding deglex order on the monomials in $x$ and $y$, and the reduction system
$$\mathcal{R}=\{x^2,\, y^2-xy,\, yxy\}.$$
We obtain that $S(A)=\bigoplus_{p\geq 0} A\otimes k.S_p \otimes A$, where $k.S_p$ denotes the $k$-vector space generated by the set $S_p$. 
Explicitely, $S_0=\{1\}$, $S_1=\{x,\,y\}$ and $S_2=\{x^2,\,y^2,\,yxy\}$ -- denoted by $S$ in~\cite{cs:reso}. For each $p\geq 3$, $S_p$ is the 
set  of the ($p-1$)-ambiguities defined by $S_2$. The $p$-ambiguities are the monomials obtained as minimal proper superpositions of $p$ elements 
of $S_2$. For example, $S_3=\{x^3,\, y^3,\, yxy^2,\, y^2xy,\, yxyxy\}$ and
$$S_4=\{x^4,\, y^4,\, yxy^3,\, y^3xy,\, y^2xy^2,\, yxy^2xy,\, y^2xyxy,\, yxyxy^2,\, yxyxyxy\}.$$
The differential $d$ is defined in $S_2$ by $d(1\otimes x^2 \otimes 1)=x\otimes x \otimes 1 + 1\otimes x \otimes x$, and
$$d(1\otimes y^2 \otimes 1)=y \otimes y \otimes 1 + 1 \otimes y \otimes y - x \otimes y \otimes 1 - 1 \otimes x \otimes y,$$
$$d(1 \otimes yxy \otimes 1)=yx \otimes y \otimes 1 + y\otimes x \otimes y + 1 \otimes y \otimes xy.$$
For any $p\geq 3$, $x^p$ belongs to $S_p$ and $d(1\otimes x^p \otimes 1)=x\otimes x^{p-1} \otimes 1 + (-1)^p 1\otimes x^{p-1} \otimes x$. 
Therefore, the morphism of graded $A$-bimodules $\chi: K(A) \rightarrow S(A)$ defined by the identity map on the generators of all the spaces 
$W_p$, except $y^2-xy$ which is sent to $y^2$, is a morphism of complexes, allowing us to view $K(A)$ as a subcomplex of the 
resolution $S(A)$. The proof of the following is omitted; 
it lies on rather long computations.
\Bpo Let $A=k\langle x,y \rangle /\langle x^2, y^2-xy\rangle$ be the algebra considered in this section.
\begin{itemize}
\item $H(\tilde{\chi})_2: HK_2(A) \rightarrow HH_2(A)$ is an isomorphism.
\item $HH_3(A)$ is $3$-dimensional, generated by the classes of $1\otimes x^3$, $y\otimes y^3 + 1\otimes yxy^2$ and $xy\otimes y^3 + y\otimes yxy^2$.
 Moreover, $H(\tilde{\chi})_3: HK_3(A) \rightarrow HH_3(A)$ sends $[1\otimes x^3]$ to itself. In particular, $H(\tilde{\chi})_3$ is injective and 
 not surjective.
\item $HH^2(A)$ is $2$-dimensional, generated by the classes of $1\otimes x^{\ast 2}+ y \otimes y^{\ast}x^{\ast}y^{\ast}$ and $x\otimes  x^{\ast 2} - y\otimes  y^{\ast 2}$. 
Moreover, $H(\chi^{\ast})_2: HH^2(A) \rightarrow HK^2(A)$ sends the first one to the class of $1\otimes x^{\ast 2}$, and the second one to the class 
of $-\frac{1}{2} y\otimes y^{\ast 2}$. In particular, $H(\chi^{\ast})_2$ is injective and not surjective.
\item  $HH^3(A)$ is $1$-dimensional, generated by the class of $xyx \otimes y^{\ast}x^{\ast}y^{\ast 2} + xyx \otimes y^{\ast 2}x^{\ast}y^{\ast}$. 
Moreover, $H(\chi^{\ast})_3=0$.
\end{itemize}
\Epo

\vspace{0.5 cm} \textsf{Roland Berger: Univ Lyon, UJM-Saint-\'Etienne, CNRS UMR 5208, Institut Camille Jordan, F-42023, Saint-\'Etienne, France}

\emph{roland.berger@univ-st-etienne.fr}\\

\textsf{Thierry Lambre: Laboratoire de Math\'ematiques Blaise Pascal, UMR 6620 CNRS \& UCA, Campus universitaire des C\'ezeaux, 3 place Vasarely, TSA 60026, CS 60026, 63178 Aubi\`ere Cedex, France}

\emph{thierry.lambre@uca.fr}\\

\textsf{Andrea Solotar: IMAS and Dto de Matem\'{a}tica, Facultad de Ciencias Exactas y Naturales,
Universidad de Buenos Aires, Ciudad Universitaria, Pabell\`{o}n 1,
(1428) Buenos Aires, Argentina}

\emph{asolotar@dm.uba.ar}

\end{document}